\documentclass[11pt, reqno]{amsart}

\usepackage[margin=1in]{geometry}
\usepackage{latexsym}
\usepackage{graphicx}
\usepackage{xypic}
\usepackage{hyperref}
\usepackage{subcaption} 
\usepackage{tikz-cd}
\usepackage{lineno}
\usepackage{enumitem}
\hypersetup{%
    colorlinks,
    linkcolor={red!50!black},
    citecolor={green!50!black},
    urlcolor={blue!50!black}
}
\usepackage{algorithm2e}
\SetKwComment{Comment}{/* }{ */}
\RestyleAlgo{ruled}
\SetKwInput{KwOut}{Output}

\newtheorem{theorem}{Theorem}

\newtheorem{claim}[theorem]{Claim}

\newtheorem{conjecture}[theorem]{Conjecture}
\newtheorem{corollary}[theorem]{Corollary}

\newtheorem{definition}[theorem]{Definition}

\newtheorem{lemma}[theorem]{Lemma}

\newtheorem{problem}[theorem]{Problem}
\newtheorem*{problem-n}{Problem}
\newtheorem{proposition}[theorem]{Proposition}
\newtheorem{prop}[theorem]{Proposition}

\numberwithin{equation}{section}
\numberwithin{theorem}{section}
\numberwithin{case}{section}

\numberwithin{subcase}{case}

\def\E{\mathcal{E}}

\def \a{\alpha}
\def\eps{\varepsilon}

\def \bi{\binom}
\def \cP{\mathcal{P}}

\def \bfi{\mathbf{i}}
\def \bfu{\mathbf{u}}
\def \bfv{\mathbf{v}}
\def \bfw{\mathbf{w}}

\def \PM{\textbf{PM}}
\def\a{\alpha}
\def\b{\beta}

\def\c{c^*_{k,\ell}}
\def\COMMENT#1{}
\let\COMMENT=\footnote

\title{On the Keevash-Knox-Mycroft Conjecture}
\author{Luyining Gan}
\address{School of Mathematical Sciences\\ Beijing University of Posts and Telecommunications\\ Beijing\\ China}
\email{elainegan@bupt.edu.cn}
\thanks{LG is partially supported by National Natural Science Foundation of China (12401446). JH is partially supported by National Natural Science Foundation of China (12371341).}

\author{Jie Han}
\address{School of Mathematics and Statistics\\ Beijing Institute of Technology\\ Beijing\\ China}
\email{han.jie@bit.edu.cn}


\begin{document}

\maketitle

\begin{abstract}
Given $1\le \ell <k$ and $\delta\ge0$, let $\textbf{PM}(k,\ell,\delta)$ be the decision problem for the existence of perfect matchings in $n$-vertex $k$-uniform hypergraphs with minimum $\ell$-degree at least $\delta\binom{n-\ell}{k-\ell}$.
For $k\ge 3$, $\textbf{PM}(k,\ell,0)$ was one of the first NP-complete problems by Karp.
Keevash, Knox and Mycroft conjectured that $\PM(k, \ell, \delta)$ is in P for every $\delta > 1-(1-1/k)^{k-\ell}$ and verified the case $\ell=k-1$.

In this paper we show that this problem can be reduced to the study of the minimum $\ell$-degree condition forcing the existence of fractional perfect matchings. 
Together with existing results on fractional perfect matchings, this solves the conjecture of Keevash, Knox and Mycroft for $\ell\ge 0.4k$.
Moreover, we also supply an algorithm that outputs a perfect matching, provided that one exists.
%
%
\end{abstract}


\section{Introduction}

Matchings have attracted a great deal of attention from both mathematicians and theoretical computer scientists, and are arguably the most natural extension of graph objects to hypergraphs.
However, the study of hypergraph matching problems is still a challenging task.
One particular reason for this is that finding maximal matchings in $k$-uniform hypergraphs for $k\ge 3$ 
is famously NP-complete~\cite{Karp}, in contrast to the tractability in the graph case (Edmonds' blossom algorithm~\cite{Edmonds}).

Hypergraph matchings also find exciting applications in other fields, e.g.~the Existence Conjecture of Block Designs~\cite{Keevash_design, GKLO_design}, Ryser's Conjecture on Latin Squares and Samuels' Conjecture in Probability Theory.
For applications on practical problems, one prime example is that Asadpour, Feige and Saberi~\cite{AUA} used hypergraph perfect matchings to study the \emph{Santa Claus problem}.

In this paper we continue the study of the decision problem of perfect matchings in dense hypergraphs, initiated by Karpi\'nski, Ruci\'nski and Szyma\'nska~\cite{KRS10}.
Given $k\ge 2$, a \emph{$k$-uniform hypergraph} (or $k$-graph) $H$ consists of a vertex set $V(H)$ and an edge set $E(H)$, where each edge in $E(H)$ is a set of $k$ vertices of $H$.
A subset $M\subseteq E(H)$ is a {\it{matching}} if every two edges from $M$ are vertex-disjoint.
A matching in $H$ is called \emph{perfect} if it covers all vertices of $H$.
Given a $k$-graph $H$ with an $\ell$-element vertex set $S$ (where $0\le \ell \le k-1$) we define $\deg_H(S)$ to be the number of edges containing $S$.
The {\it{minimum $\ell$-degree}} $\delta_\ell(H)$ of $H$ is the minimum of $\deg_H(S)$ over all $\ell$-element sets of vertices in $H$.

The following decision problem was raised by Keevash, Knox and Mycroft~\cite{KKM2015}, generalizing a problem of Karpi\'nski, Ruci\'nski and Szyma\'nska~\cite{KRS10} for the case $\ell=k-1$.
\begin{problem}
Given integers $\ell<k$ and $\delta\in [0,1]$, denote by $\PM(k, \ell, \delta)$ the problem of deciding whether there is a perfect matching in a given $k$-graph $H$ on $n\in k\mathbb N$ vertices with $\delta_\ell(H) \ge \delta \binom{n-\ell}{k-\ell}$.
What is the computational complexity of $\PM(k, \ell, \delta)$?
\end{problem}

The motivating fact is that for $k\ge 3$, $\PM(k, \ell, 0)$ is equivalent to the problem for general $k$-graphs, so is NP-complete; on the other hand $\PM(k, \ell, \delta)$ is trivially in P when $\delta$ is large (e.g., when $\delta>1-1/k$ by a result from~\cite{HPS}) because all such $k$-graphs contain perfect matchings.
Therefore, it is natural to ask for the point where the behavior changes.
A reduction of Szyma\'{n}ska~\cite{Szy2013} showed that $\PM(k, \ell, \delta)$ is NP-complete for $k\ge 3$ and $\delta < 1-(1-1/k)^{k-\ell}$.
In a breakthrough paper, Keevash, Knox and Mycroft~\cite{KKM2015} conjectured that $1-(1-1/k)^{k-\ell}$ is the turning point and verified the case $\ell=k-1$.

\begin{conjecture}[Keevash, Knox and Mycroft~\cite{KKM2015}]\label{conj_main}
For $1\le \ell<k$, $\PM(k, \ell, \delta)$ is in P for every $\delta >1-(1-1/k)^{k-\ell}$.
\end{conjecture}

Recently, Han and Treglown~\cite{HT2020} showed that the conjecture holds for $0.5k\le \ell\le (1+\ln(2/3))k\approx 0.59k$.
In this paper we verify Conjecture~\ref{conj_main} for all $\ell \ge 0.4k$.
In fact, our main result reduces the conjecture to the study of the minimum-degree-type threshold for the existence of a perfect fractional matching in $k$-graphs.
To illustrate this, we introduce the following definitions.

Given a $k$-graph $H=(V, E)$, a \emph{fractional matching} in $H$ is a function $\omega: E \to [0,1]$ such that for each $v\in V$ we have that $\sum_{e \ni v} w(e)\leq1$. Then $\sum _{e \in E} w(e)$ is the \emph{size} of $w$. If the size of $w$ in $H$ is $n/k$ then we say that $w$ is a \emph{perfect fractional matching}.
Given $k,\ell \in \mathbb N$ such that $\ell \leq k-1$, define $c^*_{k,\ell}$ to be the smallest number $c$ such that every $k$-graph $H$ on $n$ vertices with $\delta _{\ell} (H) \geq (c +o(1)) \binom{n-\ell}{k-\ell}$ contains a perfect fractional matching.
The following is our main result.

\begin{theorem}\label{thm_main}
Suppose  $k,\ell \in \mathbb N$ such that $1\le \ell \leq k-1$.
Then for any $\delta \in (c^*_{k, \ell},1]$, ${\bf PM}(k, \ell, \delta)$ is in $P$. 
That is, for any $\delta \in (c^*_{k, \ell},1]$, there exists a constant $c=c(k)$ such that there is an algorithm with running time $O(n^c)$
which given any $n$-vertex $k$-graph $H$ with $\delta_\ell(H) \ge \delta \binom{n-\ell}{k-\ell}$, either outputs a perfect matching of $H$, or a certificate that none exists.
\end{theorem}

In fact, in~\cite{HT2020}  a similar result was proved for $\delta\in (\delta^*, 1]$ where $\delta^*=\max\{c^*_{k, \ell},1/3\}$.
Comparing with their result, Theorem~\ref{thm_main} drops the extra 1/3 and thus extends the result to large values of $\ell$, namely, to $\ell > (1+\ln(2/3))k$.
In the conference version of this paper~\cite{GanHan-conf}, we prove the main result, Theorem~\ref{thm_main}, only for the decision problem. In the current version we also provide a polynomial-time algorithm that actually finds the perfect matching given that one exists. This improvement includes derandomising several parts of the arguments used in~\cite{GanHan-conf} using a result of Garbe and Mycroft~\cite{GaMy}, and in particular an algorithmic proof of Theorem~\ref{thm:alm_mat} which originally we just quoted from~\cite{CGHW21}.

For the parameter $c^*_{k,\ell}$, Alon, Frankl, Huang, R\"odl, Ruci\'nski, and Sudakov~\cite{AFHRRS} in 2012 made the following conjecture.
\begin{conjecture}\cite[Conjecture 1.1]{AFHRRS}
\label{conj11}
For all $\ell, k \in \mathbb N$, $c^*_{k,\ell}=1-(1-1/k)^{k-\ell}$.
\end{conjecture}
They~\cite{AFHRRS} verified the case $k-\ell\le 4$.
The conjecture was further validated by K\"uhn, Osthus and Townsend~\cite[Theorem 1.7]{KOT2014} for $\ell \geq k/2$ and by Han~\cite[Theorem 1.5]{Han2016} for $\ell = (k-1)/2$.
In a recent work, Frankl and Kupavskii~\cite[Corollary 24]{FrKu18} verified this conjecture for $\ell\ge 0.4k$.
Unfortunately, despite the efforts from experts in the field, Conjecture~\ref{conj11} is still open and appears to be very challenging for small values of $\ell$.
In fact, Conjecture~\ref{conj11} is also closely related to an old conjecture of Erd\H{o}s on the size of the largest matching in hypergraphs (in particular, results of~\cite{FrKu18, Han2016} are corollaries of the corresponding progress on the conjecture of Erd\H{o}s).

Combining Theorem~\ref{thm_main} with the current status on $c^*_{k,\ell}$ we get the following corollary.

\begin{corollary}
Conjecture~\ref{conj_main} holds for $\ell \ge 0.4k$ and for $k-\ell\le 4$.
\end{corollary}

Thus, by Theorem~\ref{thm_main}, Conjecture~\ref{conj_main} holds for all cases when $c^*_{k,\ell}=1-(1-1/k)^{k-\ell}$, that is, whenever Conjecture~\ref{conj11} holds.
Indeed, it has been proved in~\cite{CGHW21} that if $\delta > c^*_{k,\ell}$, then the $k$-graph contains a matching that covers all but exactly $k$ vertices (see Theorem~\ref{thm:alm_mat}).
Given this, our result can be viewed as the efficient detection of a certain class of \emph{divisibility constructions} that prevent the existence of perfect matchings.
As a consequence, we reduce the decision problem to an extremal problem on the existence of a perfect fractional matching, which can be recognized as a resolution on the ``computational complexity'' aspect of this problem.

We now give an overview of minimum-degree-type conditions as well as the relevant divisibility constructions.

\subsection{Minimum degree conditions and divisibility barriers}

The minimum degree conditions forcing a perfect matching have been studied extensively over the last two decades.
Focusing on the asymptotic thresholds, all known results support the following conjecture raised by H\`an, Person and Schacht~\cite{HPS}.
Note that this corresponds to the case when the decision problem is trivially in P (a trivial algorithm that always outputs yes).

\begin{conjecture}[H\`an--Person-Schacht,~\cite{HPS}]
\label{conj_PM}
Given $1\le \ell <k$ and $n\in \mathbb N $ such that $k\mid n$, if a $k$-graph $H$ on $n$ vertices satisfies $\delta_\ell(H)\ge (\max\{1/2, 1-(1-1/k)^{k-\ell}\}+o(1))\binom{n-\ell}{k-\ell}$, then $H$ contains a perfect matching.
\end{conjecture}

The conjecture has attracted a great deal of attention and so far has been verified for $\ell \ge 3k/8$ by Frankl and Kupavskii~\cite{FrKu18} and a handful of pairs of small values of $(k,\ell)$.
Note that this conjecture is slightly weaker than our problem, as e.g. for certain values of $\ell$, it \emph{suffices} to show that $c_{k,\ell}^*\le 1/2$, rather than determining the precise value of $c_{k,\ell}^*$ (and this is the reason that the known record on Conjecture~\ref{conj_PM} by~\cite{FrKu18} is slightly wider than that for the conjecture on $c_{k,\ell}^*$).

In fact under the assumption $\delta_\ell(H)\ge (c_{k,\ell}^*+o(1))\binom{n-\ell}{k-\ell}$, Chang, Ge, Han and Wang~\cite{CGHW21} recently proved that one can find a matching in $H$ of size $n/k-1$ (see Theorem~\ref{thm:alm_mat}).
However, such $H$ may or may not have a perfect matching, and, prior to this work, it is not clear how to characterize these two types of $k$-graphs.
To understand this, what is interesting to our problem is the divisibility constructions that achieve the bound $1/2$ in the above conjecture.
Consider an $n$-vertex set $V$ with a bipartition $X\cup Y$, where $X$ and $Y$ have almost equal size subject to $|Y|$ being odd.
Now define a $k$-graph $H_0$ on $V$ with the edge set consisting of all $k$-tuples that contain an even number of vertices in $Y$.
It is not hard to see that $\delta_\ell(H_0)\approx \frac12\binom{n-\ell}{k-\ell}$ and $H_0$ has no perfect matching.
To see this, note that any matching in $H_0$ covers an even number of vertices in $Y$, so not $Y$ entirely.

One can actually construct such partitions for an arbitrary number of parts.
For certain sizes of parts, divisibility conditions similar to the parity issue in the above example prevent the existence of perfect matchings.
Thus, our result and algorithm can be viewed as efficient detection of such constructions.
Indeed, in the Keevash--Knox--Mycroft proof~\cite{KKM2015} of Conjecture~\ref{conj_main} for $\ell=k-1$, they designed efficient algorithms to exhibit a number of ($O(n^{k+1})$) such partitions and tested the divisibility (solubility) for each of them.
In contrast, we show that one can focus on one partition and prove a sufficient and necessary condition for the existence of a perfect matching solely on that partition.
This will be made clear in Section 2.

\subsection{Related work}
The decision problem for perfect matchings in dense hypergraphs was first raised by Karpi\'nski, Ruci\'nski and Szyma\'nska~\cite{KRS10} for the case $\ell=k-1$, where they formulated the problem as $\PM(k, \delta)$ which is equivalent to $\PM(k, k-1, \delta)$ in this paper.
They showed that $\PM(k, 1/2-\eps)$ is in P for some absolute $\eps>0$, thus showing that $1/2$ is not the turning point for the change of behavior, while Szyma\'nska's~\cite{Szy2013} reduction showed that $\PM(k, \delta)$ is NP-complete when $\delta<1/k$.
This left a hardness gap for $\delta\in [1/k, 1/2)$.
Significant progress was made by Keevash--Knox--Mycroft~\cite{KKM_abs, KKM2015} who showed that $\PM(k, \delta)$ is in P for $\delta>1/k$.
This hardness problem was fully settled by Han~\cite{Han14_poly} who proved that $\PM(k,1/k)$ is in P. 
Very recently, this result was strengthened by Han and Keevash~\cite{HaKe20}, who showed that the minimum $(k-1)$-degree condition can be weakened to $n/k-c$ for any constant $c>0$ and their algorithm can actually output the perfect matching, provided that one exists.

The similar decision problem for Hamilton cycles (spanning cycles) has also been studied.
First, it is well-known that it is NP-complete to determine if a (2-)graph has a Hamilton cycle.
A $k$-graph $C$ is called a tight cycle if its vertices can be listed in a cyclic order so that the edges are all consecutive $k$-tuples.
For tight Hamilton cycles in dense $k$-graphs under minimum $(k-1)$-degree condition, it is known that the interval on the minimum $(k-1)$-degree condition for \emph{being non-trivially in $P$} has a sublinear length, unlike the interval $[n/k, n/2]$ for the perfect matching case. 
Indeed, it was shown by R\"odl, Ruci\'nski and Szemer\'edi~\cite{RRS08} that for an $n$-vertex $k$-graph $H$, if $\delta_{k-1}(H)\ge (1/2+o(1))n$ then $H$ contains a tight Hamilton cycle, i.e., the decision problem is trivially in P (the $o(n)$ term is eliminated for $k=3$~\cite{RRS11}); on the other hand, Garbe and Mycroft~\cite{GaMy} showed that there exists a constant $C$ such that if $\delta_{k-1}(H)\ge n/2-C$, then the decision problem of tight Hamilton cycles is NP-complete.
However, such an interval of linear length is shown to exist for looser cycles~\cite{GaMy}.

Han and Treglown~\cite{HT2020} considered the similar decision problem for $F$-factors\footnote{Given $k$-graphs $F$ and $H$, an $F$-factor in $H$ is a set of vertex-disjoint copies of $F$ whose union covers $V(H)$.} in graphs and $k$-graphs.
In particular, they determined the turning point for the $F$-factor problem for graphs and thus gave a negative answer to a question of Yuster~\cite{Yuster07}.


\section{A partition lemma and a structural theorem}

To prove Theorem~\ref{thm_main}, we shall establish a structural theorem (Theorem~\ref{thm_structural}) for perfect matchings.
Namely, 
we exhibit a sufficient and necessary condition for the existence of perfect matchings in $k$-graphs of large enough minimum $\ell$-degree, which, in addition, can be checked in polynomial time.
The heart of the proof of the structural theorem is the lattice-based absorption method developed by Han~\cite{Han14_poly}, which features a vertex partition of the given $k$-graph (Lemma~\ref{lem:partition}).


For the absorption property we need for building perfect matchings, we will make use of the key definition of \emph{reachability}, a concept that was introduced by Lo and Markstr\"om~\cite{LM1}.

\subsection{Reachability}\label{sub23}
Let $H$ be an $n$-vertex $k$-graph.
For $i\in \mathbb N$ and $\beta\in (0,1)$, we say that two vertices $u$ and $v$ in $V(H)$ are \emph{$(\beta, i)$-reachable in $H$} if there are at least $\beta n^{i k-1}$ $(i k-1)$-sets $S$ such that both $H[S\cup \{u\}]$ and $H[S\cup \{v\}]$ have perfect matchings. 
We refer to such a set $S$ as a \emph{reachable $(ik-1)$-set  for $u$ and $v$}.
We say a vertex set $U\subseteq V(H)$ is \emph{$(\beta, i)$-closed in $H$} if any two vertices $u,v\in U$ are $(\beta, i)$-reachable in $H$.
Given any $v \in V(H)$, define $\tilde{N}_{\b, i}(v,H)$ to be the set of vertices in $V(H)$ that are $(\beta,i)$-reachable to $v$ in $H$.

\subsection{Index vector and robust vector}
Given an $n$-vertex $k$-graph $H$ and integer $r\ge 0$, let $\cP =\{V_0, V_1, \dots, V_r\}$ be a partition of $V(H)$ into disjoint vertex sets, namely, $\bigcup_{0\le i\le r}V_i=V(H)$.
In this paper, every partition has an implicit ordering of its parts.

Next we introduce the index vectors and edge-lattices.
Given a $k$-graph $H$ and a partition $\cP=\{V_0, V_1, \dots, V_s, V_{s+1}, \dots, V_r\}$ of $V(H)$, the \emph{index vector} $\bfi_{\cP}(e)\in \mathbb{Z}^{r}$ of an edge $e\in E(H)$ with respect to $\cP$ is the vector whose coordinates are the sizes of the intersections of $e$ with each part of $\cP$ \emph{except $V_0$}, namely, $\bfi_{\cP}(e)|_i=|e\cap V_i|$ for $i\in [r]$, where $\bfv|_i$ is defined as the $i$th coordinate of $\bfv$. 
For any $\bfv=\{v_1,\dots, v_r\}\in \mathbb{Z}^{r}$, let $|\bfv|:=\sum_{i=1}^r v_i$.
Here we say that $\bfv\in \mathbb{Z}^{r}$ is a \emph{$k$-vector} if it has non-negative coordinates and $|\bfv| = k$.
In previous work, for $\mu>0$, the set of $\mu$-robust vectors (denoted by $I_{\cP}^\mu(H)$) is defined as the vectors $\bfi\in \mathbb{Z}^{r}$ such that $H$ contains at least $\mu n^k$ edges whose index vectors are equal to $\bfi$.
In this paper we need a more detailed description of robust vectors -- where we need to distinguish the roles of two different groups of $V_i$.

\begin{definition}[$\mu$-robust vectors]
\label{def:robust}
For $r>s\ge 0$,
{let $\cP=\{V_0, V_1, \dots, V_s, V_{s+1}, \dots, V_r\}$ be a partition of $V(H)$}.
Given $\mu>0$, define {$I_{\cP}^{s,\mu}(H):=I_{\cP,1}^{s,\mu}(H)\cup I_{\cP,2}^{s,\mu}(H)$} as the union of the following two sets:
\begin{enumerate}[label=$(\arabic*)$]
	\item the set $I_{\cP,1}^{s,\mu}(H)$ consists of all $k$-vectors $\bfi\in \mathbb{Z}^{r}$ such that $\bfi|_j=0$ for $j\in [s]$, $\sum^r_{j=s+1} \bfi|_j =k$ and $H$ contains at least $\mu n^k$ edges $e$ with $\bfi_\cP(e) = \bfi$;
	\item the set $I_{\cP,2}^{s,\mu}(H)$ consists of all $k$-vectors $\bfi\in \mathbb{Z}^r$ such that $\bfi|_i=1$ for exactly one $i\in [s]$, $\sum^r_{j=s+1} \bfi|_j =k-1$ and every vertex $v\in V_i$ is in at least $\mu n^{k-1}$ edges $e$ with $\bfi_\cP(e) = \bfi$.
\end{enumerate}
\end{definition}

The novel ingredient of this definition is the assumption (2), which helps us to classify the vertices that are 1-reachable to few vertices of $H$ (such vertices come from Lemma~\ref{S-closed} and do not exist when $\ell=k-1$).

Now we are ready to state our partition lemma, which outputs a refined partition compared to the partition lemmas in~\cite{Han14_poly, HT2020}.
Throughout the paper, 
we write $\alpha \ll \beta \ll \gamma$ to mean that 
it is possible to choose the positive constants
$\alpha, \beta, \gamma$ from right to left. More
precisely, there are increasing functions $f$ and $g$ such that, given
$\gamma$, whenever we choose some $\beta \leq f(\gamma)$ and $\alpha \leq g(\beta)$, the subsequent statement holds. 
Hierarchies of other lengths are defined analogously.

\begin{lemma}\label{lem:partition}
Given integers $k\ge 3$, $C>0$ and a real $\delta>0$, suppose we have $1/n_0 \ll \mu \ll \beta  \ll \delta' \ll \delta, 1/k, 1/C$.
Given an $n$-vertex $k$-graph $H$ with $n\ge n_0$ and $\delta_\ell(H) \ge\delta \binom{n-\ell}{k-\ell}$, there is a partition $\cP$ of $V(H)$ as
	\[
	\cP = \{V_0, V_1, \dots, V_{s}, V_{s+1}, \dots, V_r\}
	\]
such that with $c:=\lfloor 1/\delta \rfloor$
\begin{enumerate}[label=$(\arabic*)$]
\item $s\le 2^{\binom{c+k-2}{k-1}}$ and $r-s\le c$, \label{item:221}
\item $|V_0|\le k^{2^{\binom{2k-2}{k-1}}} \left(k(2k+1)^k + \binom{2k-2}{k-1} C \right)$ and $|\bigcup_{0\le i\le s}V_i|\le c\delta' n$,
\item for $1\le i\le s$, $|V_i|\ge (k-1)|V_0|+k(2k+1)^k + \binom{2k-2}{k-1} C$,  \label{item:223}
\item for $1\le i\le s$, there exists $\bfi\in I_{\cP,2}^{s,\mu}(H)$ such that $\bfi|_i=1$,  \label{item:4}
\item for $s+1\le i\le r$, $|V_i|\ge \delta' n/2$ and $V_i$ is $(\beta, 2^{c})$-closed in $H[\bigcup_{s+1\le i\le r}V_i]$.  \label{item:225}
\end{enumerate}
In particular, such a partition $\cP$ of $H$ can be found in time $O(n^{2^{c}k+1})$.
\end{lemma}

Let us briefly look at the properties (1)--(5) in Lemma~\ref{lem:partition}.
Clearly, (1) upper bounds $s$ and $r$, the number of parts of $\mathcal P$.
Item (2) upper bounds $|V_0|$ by an absolute constant, which allows us to search for the matching covering $V_0$ by brute force.
Item (4) guarantees that we can proceed with absorption for vertices in each $V_i$, $i\in [s]$.
By (3), the sizes $|V_1|, \dots, |V_s|$ are not too small, which provides enough room for storing robust edges for all robust edge-vectors in $I_{\mathcal P,2}^{s, \mu}(H)$ in the proof of Theorem~\ref{thm_structural}.
At last, (5) is required for absorption of the vertices in $\bigcup_{s+1\le i\le r}V_i$, as in previous proofs.



\subsection{Lattices, solubility and the structural theorem}\label{sub:lattices}
Keevash, Knox and Mycroft~\cite{KKM2015} introduced the following notions, which help us to transfer the divisibility problem to an algebraic setting as follows.
Let $0\le s < r$ and let $\cP = \{V_0, V_1, \dots, V_r\}$ be a partition of $V(H)$ for a $k$-graph $H$.


We define $L_{\cP}^{s, \mu}(H)$ as the \emph{lattice} (additive subgroup) in $\mathbb{Z}^r$ generated by $I_{\cP}^{s,\mu}(H)$. We write $L_{\max}^r$ for the lattice generated by all $k$-vectors, that is, $L_{\max}^r:=\{ \bfv\in \mathbb{Z}^r: |\bfv|\in k\mathbb Z \}$.

Suppose $L\subseteq L_{\max}^{r}$ is a lattice in $\mathbb{Z}^{r}$.
The \emph{coset group} of $(\cP, L)$ is $Q=Q(\cP, L):=L_{\max}^{r}/L$.
For any $\bfi\in L_{\max}^{r}$, the \emph{residue} of $\bfi$ in $Q$ is $R_Q(\bfi):=\bfi+L$. For any $A\subseteq V(H)$ of size divisible by $k$, the \emph{residue} of $A$ in $Q$ is $R_Q(A):=R_Q(\bfi_{\cP}(A))$.

Let $q \in \mathbb N$.
A (possibly empty) matching $M$ in $H$ of size at most $q$ is a \emph{$q$-solution for} $(\cP, L_{\cP}^{s, \mu}(H))$ (in $H$) if $\bfi_{\cP}(V(H)\setminus V(M))\in L_{\cP}^{s, \mu}(H)$; we say that $(\cP, L_{\cP}^{s, \mu}(H))$ is \emph{$q$-soluble} if it has a $q$-solution. 
We also need a strengthening of this definition as follows.
Given a set $U\subseteq V(H)$, we define that $(\cP, L_{\cP}^{s, \mu}(H))$ is \emph{$(U, q)$-soluble} if there is a matching $M$ in $H$ such that $M$ covers $U$ and $M$ is a $(|U|+q)$-solution.

In our proof, we shall pick a suitable $\mu>0$ and let $q$ be an upper bound of the order of the coset group $Q=L_{\max}^{r}/L_{\cP}^{s, \mu}(H)$ and $U$ be the part $V_0$.
Then we show that $H$ has a perfect matching if and only if $(\cP, L_{\cP}^{s, \mu}(H))$ is $(V_0, q)$-soluble.




\begin{theorem}[Structural Theorem]\label{thm_structural}
Let $k,\ell, q \in \mathbb N$ where $\ell\leq k-1$ and let $\gamma >0$ be given.
There exist $n_0, C:=C(k,q) \in \mathbb N$ and $\delta', \b,\mu >0$ such that
\begin{equation*}
    	1/n_0 \ll \b, \mu \ll \delta' \ll \gamma, c_{k,\ell}^*, 1/q, 1/C, 1/k.
\end{equation*}
Let $H$ be an $n$-vertex $k$-graph with $\delta_\ell(H)\ge (c_{k,\ell}^*+\gamma)\binom{n-\ell}{k-\ell}$, where $n \geq n_0$ and  $k$ divides $n$.
Suppose $\cP$ is a partition of $V(H)$ satisfying Lemma~\ref{lem:partition} $(1)$--$(5)$ with $\delta=c_{k,\ell}^*$.
Moreover, suppose $|Q(\cP, L^{s, \mu}_{\cP}(H))|\leq q$.
Then $H$ contains a perfect matching if and only if $(\cP, L^{s, \mu}_{\cP}(H))$ is $(V_0, q)$-soluble.
Furthermore, there is an algorithm with running time $O(n^{2^{k+2}k^2})$ that finds a perfect matching in $H$ if $(\cP, L^{s, \mu}_{\cP}(H))$ is $(V_0, q)$-soluble.
\end{theorem}


\section{Highlights of the proof: a comparison with the Han--Treglown proof} \label{sec:comparison}

The basic idea for establishing the structural theorem is to distinguish the roles of \emph{robust} and \emph{non-robust} edges: to avoid the divisibility barriers, we may have to use edges with certain (combinations of) index vectors.
For some index vectors $\bfv$ there are many edges $e$ with $\bfi_\cP(e)=\bfv$, namely, there are many ``replacements'' even when we are forbidden from using, say, a small number of such edges.
For other index vectors $\bfv$ there are few edges $e$ with $\bfi_\cP(e)=\bfv$, so we have to be careful when using such edges. 
In fact, the algebraic setting allows us to show that one can restrict the attention to only a \emph{constant} number of such non-robust edges (using the lattice and coset group arguments), and thus this can be tested by brute force.
Then the rest of the proof follows from the lattice-based absorption argument.
Roughly speaking, it reserves a small matching which can be used to turn an almost perfect matching to a perfect matching given certain divisibility conditions on the leftover vertices.

In~\cite{HT2020} Han and Treglown proved our Theorem~\ref{thm_main} under the additional assumption that $\delta > 1/3$, which verified Conjecture~\ref{conj_main} for $0.5k\le \ell \le 0.59k$ (now together with the recent result in~\cite{FrKu18} it resolves Conjecture~\ref{conj_main} for $0.4k\le \ell \le 0.59k$).
However, this does not solve the conjecture for $\ell=k-2$, which might be considered as the easiest case after the resolution of the case $\ell=k-1$.
Below we shall first outline the proof in~\cite{HT2020}, and then explain our innovation compared with their approach and how such an improvement is achieved.

The partition lemma used in~\cite{HT2020} is Lemma~\ref{lem:P} in this paper (which we use as a building block to establish our partition).
The key issue is that when $\ell<k-1$, one can not apply Lemma~\ref{lem:P} directly to the $k$-graph $H$, as in $H$ there might be a set $W$ of vertices $v$ which are not reachable to many vertices, namely, $\tilde{N}_{\b, i}(v,H)$ is small for any proper choice of $\b>0$ and $i\in \mathbb N$.
However, it is straightforward to show that $|W|=o(n)$, and (after some work) we can apply Lemma~\ref{lem:P} with $S=V(H)\setminus W$ and get a partition of $V(H)\setminus W$.
Now we face the following challenge. 
\begin{problem-n}
Suppose $|W|=o(n)$. How do we find a matching $M$ covering $W$ so that $H-V(M)$ has a perfect matching (or conclude that none exists)?
\end{problem-n}
The problem is trivial if $|W|$ is a constant, for which we can do brute force search for a matching $M$ of constant size, which involves $O(n^{|W|})$ possibilities; when $|W|$ is large, it is hopeless without further assumptions.

Furthermore, it was not clear how to deal with the vertices of $W$ by absorption, as $|W|$ might be smaller than the threshold for $\mu$-robustness but still a small linear size, i.e., $\eps n \le |W| < \mu n$, so that every vector touching $W$ will not be recorded as a $\mu$-robust vector.
The proof in~\cite{HT2020} avoided the ``decision'' part of the problem by assuming $\delta > 1/3$, so that when $W$ is non-empty $V(H)\setminus W$ is closed, in which case $H$ always contains a perfect matching (indeed it is shown that any matching $M$ covering $W$ will be a solution).
Therefore, the problem is left open for $\delta <1/3$ (i.e., for $\ell > (1+\ln(2/3))k \approx 0.595k$).

We also note that the existence of $W$ is not a problem in the existential results in the literature.
In previous work on sufficient (minimum-degree-type) conditions for perfect matchings, those vertices can be put into a matching of small linear size, whose removal does not affect much the minimum-degree conditions, guaranteeing that the absorption can proceed after the removal of this small matching.

\medskip
Our new proof can be seen as a significant refinement of the previous approach, where we strengthen our control on both the partition and the robust vectors.
As mentioned earlier, our new proof features a finer partition lemma (Lemma~\ref{lem:partition}) than previous ones, where we classify vertices in $W$ as well.
More precisely, we first partition $S:=V(H)\setminus W$, the set of vertices which are 1-reachable to $\Omega(n)$ other vertices, by Lemma~\ref{lem:P} and denote the partition by $\cP_1=\{W_1, \dots, W_d\}$.
Then we classify vertices of $W$ according to their edge distributions in $\cP_1$, that is, we obtain a partition of $W$ by collecting vertices with common robust edge vectors together, so that the partition satisfies Definition~\ref{def:robust} (2).
Next we put the clusters that are too small (smaller than a certain constant) to a trash set $V_0$ in a recursive manner.
This results in a trash set $V_0$ of constant order, and because we have no control on $V_0$ at all, we will check how to match $V_0$ by brute force in time $O(n^{|V_0|})$.
Now the relief is that all clusters that survived from this greedy process have a good (though still perhaps constant) size (Lemma~\ref{lem:partition} (3)), which is enough (and crucial) for a (refined) absorption argument to work in later proofs.
Since all the above procedures can be done in polynomial time, we get the desired polynomial-time algorithm for the decision problem $\PM(k,\ell, \delta)$.


\section{Proof of Theorem~\ref{thm_main}}

\subsection{An upper bound on the size of the coset group}
We first provide an upper bound on the size of the coset group.

\begin{proposition}
\label{prop:Q}
Suppose $1/n\ll \mu \ll \delta' \ll \gamma \ll 1/k$.
Let $H$ be a $k$-graph with $\delta_{\ell}(H) \ge (c_{k,\ell}^*+\gamma) \binom n{k-\ell}$ and $\cP=\{V_0, \dots, V_s, V_{s+1},\dots, V_{r}\}$ be a partition of $V(H)$ satisfying Lemma~\ref{lem:partition} (1)--(5) with $\delta = c_{k,\ell}^*$.
Then we have $|Q(\cP, L^{s, \mu}_{\cP}(H))|\leq (2k+1)^{r-s}$.
\end{proposition}

\begin{proof}
Write $L:=L^{s, \mu}_{\cP}(H)$ and $Q:=Q(\cP, L)$.
We first define robust $k$-vectors $\bfv_1,\dots, \bfv_r\in I^{s, \mu}_{\cP}(H)$ as follows.
For $i\in [s]$, as $\delta_1(H)\ge (1/k)\binom{n}{k-1}$, $H$ contains many edges with at least one vertex in $V_i$.
Using $|\bigcup_{0\le i\le s}V_i|\le c\delta' n$ (Lemma~\ref{lem:partition} (2)), $H$ contains at least $(1/2k)\binom nk$ edges with exactly one vertex in $V_i$ and other vertices in $V_{s+1}\cup \cdots \cup V_r$.
By the pigeonhole principle over $\binom{r-s+k-1}{k-2}$ such choices, there exists $\bfv_i\in I^{s, \mu}_{\cP,2}(H)$ with its $i$-th coordinate equal to 1 and all other of the first $s$ coordinates equal to 0.
This defines $\bfv_1,\dots, \bfv_s\in I^{s, \mu}_{\cP,2}(H)$.

For $i\in [s+1, r]$, as $\delta_{\ell}(H) \ge (c_{k,\ell}^*+\gamma) \binom n{k-\ell}$, every set of $\ell$ vertices in $V_{s+i}$ is in at least $(c_{k,\ell}^*+\gamma) \binom n{k-\ell}$ edges of $H$.
Together with $|V_{s+i}|\ge \delta' n$, $H$ contains at least $c_{k,\ell}(\delta')^\ell \binom n{k} \ge \binom{r+\ell}{\ell}\mu \binom n{k}$ edges each of which contains at least $\ell$ vertices of $V_{s+i}$.
Then as $|\bigcup_{0\le i\le s}V_i|\le c\delta' n$, by the pigeonhole principle, we obtain a $k$-vector $\bfv_i\in I^{s, \mu}_{\cP,1}(H)$ such that $\bfv_i$ has $(s+i)$-th coordinate at least $\ell$ and all its first $s$ coordinates are equal to 0.
This defines $\bfv_{s+1},\dots, \bfv_r\in I^{s, \mu}_{\cP,1}(H)$.

Note that it suffices to show that every coset contains a vector $\bfv'$ whose all coordinates have absolute values at most $k-1$ and the first $s$ coordinates equal to 0, which implies the proposition.
Given any vector $\bfv$ in any coset, we first add appropriate multiples of $\bfv_i, i\in [s]$ to $\bfv$ so that the resulting vector $\bfv^*$ have its first $s$ coordinates equal to 0 and $\bfv^*\in \bfv+L$.
Then it remains to consider the last $r-s$ coordinates.

We first note that for $\ell=1$, we have $c_{k,\ell}^* > 1/2$ and thus $r-s=1$;
and for $\ell \le 2k/3$, 
it holds that $c_{k,\ell}^* > 1/3$ and thus $r-s\le 2$.
In this case the argument below is essentially from~\cite{HT2020}.
If $r-s=1$, then $\bfv_r=k\bfu_r$ is in every coset of $L$, yielding that $|Q|=1$.
Otherwise we have $\ell\ge 2$ and $r-s=2$.
If there are two distinct $k$-vectors $\bfw_{1}, \bfw_2\in I^{s, \mu}_{\cP,1}(H)\subseteq L$, then we have $\bfw_{1}- \bfw_2 = a \bfu_{s+1} - a \bfu_r\in L$ for some nonzero $a\in [-k,k]$.
Then we can add appropriate multiples of $\bfw_1$ and  $a \bfu_{s+1} - a \bfu_r$ to $\bfv^*$ so that the resulting vector $v'\in \bfv^*+L=\bfv+L$ has the last two coordinates sum to 0 and lie in $[-|a|, |a|]$.
This yields that $|Q|\le k$.
Therefore it remains to consider the case that $I^{s, \mu}_{\cP,1}(H)$ has exactly one element, denoted by $\bfw:=a \bfu_{s+1} + (k-a) \bfu_r$.
As then we must have $\bfv_{s+1}=\bfv_r=\bfw$, we have $a\ge \ell$ and $k-a\ge \ell$.
Note that $H$ induced on $V_{s+1}\cup V_r$ has at most $\binom{|V_{s+1}|}{a}\binom{|V_{r}|}{k-a}+O(\mu n^k)$ edges.
Therefore, by averaging, for each $0\le \ell'\le \ell$, there exists an $(\ell', \ell-\ell')$-set $S_{\ell'}$ such that
\[
\deg_H(S_{\ell'})\le \binom{|V_{s+1}|}{a-\ell'}\binom{|V_{r}|}{k-a-(\ell-\ell')}+O(\mu n^{k-\ell})+O(\delta'n^{k-\ell}).
\]
Summing over $\ell'$ we have
\[
\sum_{0\le \ell'\le \ell}\deg_H(S_{\ell'})\le \binom{n}{k-\ell}+O(\delta'n^{k-\ell}).
\]
As $\ell\ge 2$, the above sum contains at least three terms and thus yields the existence of some $\ell'$ with $\deg_H(S_{\ell'})\le \frac13\binom{n}{k-\ell}+O(\delta'n^{k-\ell})$, contradicting $\delta_{\ell}(H)\ge (c_{k,\ell}^*+\gamma) \binom n{k-\ell} \ge (1/3+\gamma) \binom n{k-\ell}$ and $\delta'\ll \gamma$.


Finally we consider $\ell > 2k/3$.
Recall that $\bfv^*\in \bfv+L$ have all its first $s$ coordinates equal to 0.
Then as long as $\bfv^*$ has a coordinate in the position $s+1$ to $s+r$ with absolute value at least $k$ (say the $i$-th position), then we add $\bfv_i$ or $-\bfv_i$ to it.
This reduces the sum of the absolute values of these coordinates by at least $\ell - (k-\ell)\ge  k/3$.
We execute this procedure iteratively and obtain a vector $\bfv'\in \bfv+L$ such that all the last $r-s$ coordinates have absolute value at most $k-1$.
This yields the desired bounds on $|Q|$ and the proof is completed.
\end{proof}

Now we prove Theorem~\ref{thm_main}.
Recall that $c_{k,\ell}^*\ge c_{k,k-1}^* = 1/k$.
Then $\lfloor 1/c_{k,\ell}^*\rfloor \le k$.
Let
\[
q=q(k):= (2k+1)^k.
\]
and 
let $C:= C(k,q)$ be given by Theorem~\ref{thm_structural}.
Suppose we have constants satisfying the following hierarchy
\begin{equation*}
	1/n_0 \ll \mu \ll \b \ll \delta' \ll \gamma, c_{k,\ell}^*, 1/q, 1/C, 1/k.
\end{equation*}

Both Lemma~\ref{lem:partition} and Theorem~\ref{thm_structural} require that $n$ is larger than a constant $n_0$, and by custom $k$-graphs with less than $n_0$ vertices can be tested by brute force. 
By Lemma~\ref{lem:partition}, in time $O(n^{2^{k-1}k+1})$ we can find a partition $\cP$ satisfying Lemma~\ref{lem:partition} (1)--(5).
Because of Lemma~\ref{lem:partition} (1), we know $r-s\le k$ and obtain that $|Q(\cP, L^{s, \mu}_{\cP}(H))|\leq {\color{red}(2k+1)^{r-s}}\le q$ by Proposition~\ref{prop:Q}.
Then by Theorem~\ref{thm_structural}, it suffices to test if $(\cP, L_\cP^{s,\mu}(H))$ is $(V_0, q)$-soluble.
This can be done by testing whether any matching $M$ of size at most $|V_0|+q$ covering $V_0$ is a solution of $(\cP, L^{s, \mu}_{\cP}(H))$, in time $O(n^{k|V_0|+kq})$.
The overall time is polynomial in $n$ because $q = q(k)$
and 
\[
|V_0|\le k^{2^{\binom{2k-2}{k-1}}} \left(kq + \binom{2k-2}{k-1} C \right),
\]
where we recall that $C:= C(k, q)$ only depends on $k$.
Note that if $(\cP, L_\cP^{s,\mu}(H))$ is not $(V_0, q)$-soluble, then we output the pair $(\cP, L_\cP^{s,\mu}(H))$, which certifies that $H$ has no perfect matching (by Theorem~\ref{thm_structural}).

\medskip
\noindent\textbf{Organization.}
The rest of this paper is organized as follows.
Note that it remains to prove Lemma~\ref{lem:partition} and Theorem~\ref{thm_structural}.
We collect and prove a number of auxiliary results and give a proof of Lemma~\ref{lem:partition} in Section 5.
In Section 6, we give two absorption lemmas, which are used in the algorithmic proof of Theorem~\ref{thm:coverbut} in Section 7. 
The proof of Theorem~\ref{thm_structural} is presented in Section 8.


\section{Useful tools}
In this section we collect together some results that will be used in our proof of Theorem~\ref{thm_structural}.
When considering $\ell$-degree together with $\ell'$-degree for some $\ell'\neq \ell$, the following well-known proposition is very useful.

\begin{prop} \label{prop_deg}
Let $0\le \ell \le \ell ' < k$ and $H$ be a $k$-graph on $n$ vertices. If $\delta_{\ell '}(H)\geq x\binom{n- \ell '}{k- \ell '}$ for some $0\le x\le 1$, then $\delta_{\ell}(H)\geq x\binom{n- \ell}{k- \ell}$.
\end{prop}

This proposition is straightforward since $\delta_{\ell}(H)\geq \binom{n- \ell}{\ell'- \ell} \delta_{\ell'}(H)/{{k-\ell} \choose {\ell'-\ell}}$.



\subsection{Almost perfect matchings}
Let $k,\ell \in \mathbb N$ where $\ell\leq k-1$. 
Given $D\in \mathbb N$ and a $k$-graph $H$ on $n\in k\mathbb N$ vertices with $\delta _{\ell} (H) \geq \delta  \binom{n-\ell}{k-\ell}$, it is proved in~\cite{HT2020} that the infimum of $\delta$ such that $H$ contains a matching covering all but at most $D$ vertices is at most $\max\{1/3, c^*_{k,\ell}\}$.
We need the extra term 1/3 removed, which was very recently proved by Chang, Ge, Han and Wang~\cite{CGHW21}.

\begin{theorem} \label{thm:alm_mat}\cite{CGHW21}
Given $k, \ell\in \mathbb N$ such that $1\le \ell\le k-1$ and $\gamma >0$, there exists $n_0\in \mathbb N$ such that for $n\ge n_0$ the following holds.
Suppose $H$ is an $n$-vertex $k$-graph with $\delta_\ell(H)\ge (c_{k,\ell}^*+\gamma)\binom{n-\ell}{k-\ell}$, then $H$ contains a matching $M$ that covers all but at most $2k-\ell-1$ vertices.
In particular, when $n\in k\mathbb N$, $M$ is a perfect matching or covers all but exactly $k$ vertices.
\end{theorem}

Since the original proof of Theorem~\ref{thm:alm_mat} in~\cite{CGHW21} is not constructive, we shall give an alternative proof in Section 7 by derandomising the probabilistic proofs of several key lemmas in~\cite{CGHW21} supporting Theorem~\ref{thm:alm_mat}.

\subsection{Build a partition}
To build a partition, we need the following partition lemma from~\cite{HT2020}. 

\begin{lemma}\cite[Lemma 6.3]{HT2020}\label{lem:P}
Let $\delta'>0$, and integers $c, k\ge 2$ be given and suppose $1/n \ll \beta \ll \a \ll 1/c, \delta'$.
Assume $H$ is an $n$-vertex $k$-graph and $S\subseteq V(H)$ is such that $|\tilde{N}_{\a, 1}(v,H)\cap S| \ge \delta' n$ for any $v\in S$. Further, suppose every set of $c+1$ vertices in $S$ contains two vertices that are $(\a, 1)$-reachable in $H$. 
Then in $O(n^{2^{c-1} k+1})$ time we can find a partition $\cP$ of $S$ into $V_1,\dots, V_r$ with $r\le \min\{c, 1/\delta' \}$ such that for any $i\in [r]$, $|V_i|\ge (\delta' - \a) n$ and $V_i$ is $(\beta, 2^{c-1})$-closed in $H$.
\end{lemma}

%
%
%
%

As mentioned in Section~\ref{sec:comparison}, we cannot apply Lemma~\ref{lem:P} directly to $H$, because $H$ may contain vertices that are reachable to few other vertices. We collect them in a greedy manner in the following lemma. 
Note that a similar lemma was used in~\cite{DHSWZ21}.
\begin{lemma}\label{S-closed}
Let $\a>0$, and integers $c,k\geq 2$ be given and suppose $1/n\ll \delta' \ll \a, 1/k, 1/c$. 
Assume that $H$ is a $k$-graph on $n$ vertices satisfying that every set of $c+1$ vertices contains two vertices that are $(2\a ,1)$-reachable in $H$. 
Then in time $O(c n^{k+1})$ we can find a set of vertices $S\subseteq V(H)$ with $|S|\geq (1-c\delta' )n$ such that $|\tilde{N}_{\a, 1}(v, H[S])|\geq \delta' n$ for any $v\in S$.
\end{lemma}

We remark that in the above lemma it is important to obtain the conclusion on $\tilde{N}_{\a, 1}(v, H[S])$ rather than $\tilde{N}_{\a, 1}(v)\cap S$.
Indeed, in the latter one the reachable sets are still defined in $H$, so may contain vertices in $V(H)\setminus S$.
This is not strong enough for our later proof (see Lemma~\ref{lem:abs} and its proof).

\begin{proof}
Let $H$ be a $k$-graph on $n$ vertices satisfying the condition of Lemma~\ref{S-closed}. We greedily identify vertices with few ``reachable neighbors'' and remove the vertex together with the vertices reachable to it from $H$.
Set $V_0:=V(H)$. 
First, for every two vertices $u, v\in V(H)$, we determine if they are $(\a, 1)$-reachable in $H$, which can be done by testing if any $(k-1)$-set is a reachable set in time $O(n^{k-1})$.
Summing over all pairs of vertices, this step can be done in time $O(n^{k+1})$.
Then we check if there is a vertex $v_0\in V_0$ such that $|\tilde{N}_{\a , 1}(v_0, H)|<\delta' n$ in time $O(n^2)$.
If there exists such a vertex $v_0$, then let $A_0:=\{v_0\}\cup \tilde{N}_{\a, 1}(v_0, H)$ and let $V_1:=V_0\setminus A_0$.
Next, we check if there exists a vertex $v_1\in V_1$ such that $|\tilde{N}_{\a , 1}(v_1, H[V_1])|<\delta' n$, and if yes, then let $A_1:=\{v_1\}\cup \tilde{N}_{ \a , 1}(v_1, H[V_1])$ and let $V_2:=V_1\setminus A_1$ and repeat the procedure until no such $v_j$ exists.

Suppose when the process terminates we obtain a set of vertices $v_0,\dots, v_s$.
We claim that $s< c$ and thus $|\bigcup _{0\leq i\leq s}A_i|\leq c\delta' n$.
Indeed, otherwise consider $v_0,\dots, v_c$, the first $c+1$ of them and we shall show that every pair of them is not $(2\a , 1)$-reachable in $H$, contradicting our assumption.
Given $0\le i<j\le c$, as $v_j\notin \tilde{N}_{ \a , 1}(v_i, H[V_i])$, $v_i$ and $v_j$ have less than $\a n^{k-1}$ 1-reachable sets in $H[V_i]$.
Also, because $\delta' \ll \a, 1/c$, there are at most $c\delta' n\cdot n^{k-2}\le \a n^{k-1}$ 1-reachable sets in $H\setminus E(H[V_i])$.
These two together yield that $v_i$ and $v_j$ are not $(2\a , 1)$-reachable in $H$.

This greedy procedure needs to recompute $\tilde{N}_{ \a , 1}(v, H[V_i])$ at each time and can be done in time $O(cn^{k+1})$.
Set $S:=V(H)\setminus (\bigcup _{0\leq i\leq s}A_i )$.
We have $|S|\geq (1-c\delta' )n$ and $|\tilde{N}_{\a , 1}(v, H[S])|\geq \delta' n$ for every $v\in S$.
\end{proof}

Now we are ready to prove Lemma~\ref{lem:partition}, which establishes our vertex partition.


\subsection{Proof of Lemma~\ref{lem:partition}}

%


Choose additional constants $\alpha, \alpha'$ and $\gamma$ such that
\[
1/n_0 \ll  \mu \ll \beta  \ll \a \ll \gamma, \delta'\ll \alpha' \ll \delta, 1/k, 1/C.
\]
Assume  $n\geq n_0$ and $k$ divides $n$.
Let $H$ be an $n$-vertex $k$-graph.
Write $c:=\lfloor 1/\delta \rfloor$, then by Proposition~\ref{prop_deg} we have
\[
(c+1)\delta_1(H) \ge (c+1) \delta\bi{n-1}{k-1}>(1+\gamma)\bi{n-1}{k-1}.
\]
Thus every set of $c+1$ vertices of $V(H)$ contains two vertices that are $(2\a',1)$-reachable, as otherwise, by the inclusion-exclusion principle and $\a' \ll \gamma, \delta$
\[
n \ge (c+1)\delta_1(H) - \binom{c+1}2\cdot 2\alpha' n^{k-1} \ge (1+\gamma)\bi{n-1}{k-1} - (c+1)^2\alpha' n^{k-1} >n
\]
a contradiction.
%

By Lemma~\ref{S-closed}, we find $S\subseteq V(H)$ with $|S| \ge (1-c\delta')n$ such that $|\tilde{N}_{\a', 1}(v, H[S])|\ge \delta' n$ for every $v\in S$, in time $O(n^{k+1})$.
Let $V':= V(H)\setminus S$ and thus $|V'| \le c\delta'n$.
Now by $\a < \a'$, for every $v\in S$, we have $|\tilde{N}_{\a, 1}(v, H[S])|\ge|\tilde{N}_{\a', 1}(v, H[S])|\ge \delta' n$ and every set of $c+2$ vertices of $S$ contains two vertices that are $(2\a,1)$-reachable in $H$.
Apply Lemma~\ref{lem:P} to $H[S]$, and in time $O(n^{2^{c}k+1})$ we find a partition $\cP_1$ of $S$ into $W_{1}, \dots, W_{d}$ with $d \le c$ such that for $i\in [d]$, $|W_i|\ge (\delta'-\a)n$ and $W_i$ is $(\b, 2^{c})$-closed in $H[S]$.


Let $I^{k-1}_d$ be the set of all $(k-1)$-vectors on $\cP_1$ and note that $|I^{k-1}_d|=\binom{d+k-2}{d-1}$.
Let $\mathcal I$ be the collection of all subsets of $I^{k-1}_d$ and clearly $|\mathcal I|= 2^{|I^{k-1}_d|}=2^{\binom{d+k-2}{d-1}}$.
We classify the vertices in $V'$ by the types of the edges in which they are contained.
Indeed, for $I\in \mathcal I$, let $V_I$ be the collection of vertices $v\in V'$ such that the following two properties hold:
\begin{itemize}
\item for every $\bfi\in I$, there are at least $\mu n^{k-1}$ edges $e$ of $H$ such that $v\in e$ and $\bfi_{\cP_1}(e\setminus \{v\})=\bfi$;
\item for every $\bfi\notin I$, there are fewer than $\mu n^{k-1}$ edges $e$ of $H$ such that $v\in e$ and $\bfi_{\cP_1}(e\setminus \{v\})=\bfi$.
\end{itemize}
Clearly this defines a partition of $V'$.
Moreover, note that $V_{\emptyset}=\emptyset$ -- this is because any vertex in $V_{\emptyset}$ has vertex degree at most 
\[2^{\binom{d+k-2}{d-1}}\mu n^{k-1} + |V'|n^{k-2}\le 2^{\binom{d+k-2}{d-1}}\mu n^{k-1} + c\delta'n^{k-1}< \delta_1(H),\] 
violating the minimum degree assumption.
In particular, this implies~\ref{item:4}.
Note that this partition can be built by reading the edges for each $v\in V'$, so in time $O(n^{k})$.
Next we collect the parts that are too small and put them into a trash set $V_0$ in a recursive manner.

We first sort $V_I$, $I\in \mathcal I$ such that $|V_I|$ is increasing.
Next, starting from $V_0=\emptyset$, we \emph{recursively} check in time $O(|\mathcal I|n)$ if the next $V_I$, $I\in \mathcal I$ in the sequence satisfies that
\[
|V_{I}| < (k-1)|V_0|+ b, \text{ where } b:=k(2k+1)^k + \binom{2k-2}{k-1} C.
\]
and if yes, put all vertices of $V_I$ to $V_0$ (note here that $V_0$ is dynamic).
Because $|\mathcal I|=2^{\binom{d+k-2}{d-1}}$, straightforward computation shows that after the process we have
\[
|V_0|\le \frac{k^{|\mathcal I|}-1}{k-1} b \le k^{2^{\binom{2k-2}{k-1}}} \left(k(2k+1)^k + \binom{2k-2}{k-1} C \right).
\]

At last, in constant time we remove the empty clusters and relabel the remaining parts $V_I$ to $V_1, \dots, V_s$, and relabel the parts of $\cP_1$ as $V_{s+1}, \dots, V_{r}$.
The resulting partition satisfies all desired properties in the lemma and the overall running time is $O(n^{2^{c}k+1})$.
\section{Absorption lemmas}
In this section we prove three absorption lemmas -- the first one will be used in the proof of our main result, Theorem~\ref{thm_structural} and the other two are for in the (alternative) proof of Theorem~\ref{thm:alm_mat}.
The absorption method is by now a standard way to turn an almost spanning structure into a spanning one.
Here we use a variant called the \emph{lattice-based absorption method}, developed by Han~\cite{Han14_poly}.
Although the lemmas were originally proved by probabilistic methods, it can be derandomized by the so-called conditional expectation method.
Here we use a result of Garbe and Mycroft~\cite{GaMy}, which in turn was inspired by Karpi\'nski, Ruci\'nski and Szyma\'nska~\cite{KRS10hc}.

\begin{lemma}\cite[Proposition 4.7]{GaMy} \label{lem:derandomization}
	Fix constants $\beta>\tau>0$ and integers $m, M, N$ and $r\le N$ such that $r$ and $N$ are sufficiently large, and that $M\le (1/8)\exp(\tau^2r/(3\beta))$. Let $U$ and $W$ be disjoint sets of sizes $|U| = M$ and $|W| = N$. Let $G$ be a graph with vertex set $U\cup W$ such that $G[U]$ is empty, $G[W]$ has precisely $m$ edges, and $\deg_G(u)\ge \beta N$ for every $u\in U$. Then in time $O(N^4+MN^3)$ we can find an independent set $R\subseteq W$ in $G$ such that $(1-\nu)r\le |R|\le r$ and $|N_G(u)\cap R|\ge (\beta - \tau -\nu)r$ for all $u\in U$, where $\nu = 2mr/N^2$.
\end{lemma}

Fix an integer $i>0$. Let $H$ be a $k$-graph. For a $k$-set $S$, we says a set $T$ is an \emph{absorbing $i$-set} for $S$ if $|T|=i$ and both $H[T]$ and $H[T\cup S]$ contain perfect matchings.
Now we are ready to present our first absorption lemma, which is similar to~\cite[Lemma 3.4]{Han14_poly}.
The only difference is due to our refined definition of robust vectors $I_{\cP}^{s, \mu} (H)$. 

\begin{lemma}[Absorption Lemma] \label{lem:abs}
Suppose $k\ge 3$, $\delta >0$ and let $t:=2^{\lfloor 1/\delta \rfloor}$. Suppose that
\[
1/n_0 \ll \alpha \ll \beta, \mu \ll 1/t, 1/k.
\]
Let $H$ be an $n$-vertex $k$-graph with a partition $\cP$ of $V(H)$ satisfying Lemma~\ref{lem:partition} $(1)$-$(5)$, where $n \geq n_0$ and  $k$ divides $n$.
Let $n_1:=|\bigcup_{s+1\le i\le r} V_i|$ (where $r, s$ are from the statement of Lemma~\ref{lem:partition}).
Then in time $O(n^{4tk^2})$, we can find a family $\E_{abs}$ on $\bigcup_{s+1\le i\le r} V_i$ consisting of at most $\beta n_1$ disjoint $t k^2$-sets such that for each $A\in \E_{abs}$, $H[A]$ contains a perfect matching and every $k$-set $S\subseteq V(H)$ with $\bfi_{\cP}(S)\in I_{\cP}^{s, \mu} (H)$ has at least $\alpha n_1$ absorbing $t k^2$-sets in $\E_{abs}$. 
\end{lemma}

\begin{proof}
Roughly speaking, in the proof we first exhibit a large number of absorbing sets for each $k$-set $S$ with $\bfi_{\cP}(S)\in I_{\cP}^{s, \mu} (H)$, and then show that the desired family $\E_{abs}$ can be obtained by applying Lemma~\ref{lem:derandomization}. 
Our first task is to prove the following claim.

\begin{claim}\label{clm:abs}
Any $k$-set $S$ with $\bfi_{\cP}(S)\in I_{\cP}^{s, \mu}(H)$ has at least $\mu^{t+1} \b^{k+1} n_1^{t k^2}$ absorbing $t k^2$-sets which consist of vertices only in $\bigcup_{s+1\le i\le r}V_i$.
\end{claim}

\begin{proof}
We split the proof into two cases relating to $I_{\cP,1}^{s,\mu}(H)$ and $I_{\cP,2}^{s,\mu}(H)$ respectively.
Note that all reachable sets will be constructed with vertices in $\bigcup_{s+1\le i\le r}V_i$ only.

\noindent\textbf{Case 1.}
Suppose $\bfi \in I_{\cP,1}^{s,\mu}(H)$.
For a $k$-set $S=\{y_1,\dots, y_k\}$ with $\bfi_{\cP}(S) =\bfi$, we construct absorbing $t k^2$-sets for $S$ as follows. 
We first fix an edge $W=\{x_1, \dots, x_k\}$ in $H$ such that $\bfi_{\cP}(W)=\bfi$ and $W\cap S=\emptyset$. 
Note that we have at least $\mu n^k - k n_1^{k-1} >\frac{\mu}2 n^k$ choices for such an edge. 
Without loss of generality, we may assume that for all $i\in [k]$, $x_i, y_i$ are in the same part $V_{j}$ of $\cP$ for $j> s$.
Recall that by Lemma~\ref{lem:partition} (5), $V_{j}$ is $(\beta, t)$-closed in $H[\bigcup_{s+1\le i\le r}V_i]$.
Since $x_i$ is $(\beta, t)$-reachable to $y_i$, there are at least $\beta n_1^{t k-1}$ $(t k-1)$-sets $T_i$ such that both $H[T_i\cup \{x_i\}]$ and $H[T_i\cup \{y_i\}]$ have perfect matchings. 
We pick disjoint reachable $(t k-1)$-sets for each $x_i, y_i$, $i\in [k]$ greedily, while avoiding the previously chosen vertices. 
Since the number of previously chosen vertices is at most $t k^2+k$, we have at least $\frac{\beta}2 n_1^{t k-1}$ choices for such $(t k-1)$-sets in each step.
Note that  $W\cup T_1\cup \cdots \cup T_{k}$ is an absorbing set for $S$. 
First, it contains a perfect matching because each $T_i\cup \{x_i\}$ for $i\in [k]$ spans $t$ vertex-disjoint edges.
Second, $H[W\cup T_1\cup \cdots \cup T_{k}\cup S]$ also contains a perfect matching as each $T_i\cup \{y_i\}$ for $i\in [k]$ spans $t$ vertex-disjoint edges. 
There are at least $\frac{\mu}2 n_1^k$ choices for $W$ and at least $\frac{\beta}2 n_1^{t k-1}$ choices for each $T_i$.
Thus we find at least 
\[
\frac{\mu}2 n_1^k \times \frac{\beta^{k}}{2^k} n_1^{t k^2-k} \times \frac{1}{(tk^2)!} \geq \mu \beta^{k+1} n_1^{t k^2}
\] absorbing $t k^2$-sets for $S$.

\medskip
\noindent\textbf{Case 2.}
	Suppose $\bfi \in I_{\cP,2}^{s,\mu}(H)$.
Suppose $S=\{v_1, y_2,\dots, y_k\}$ with $\bfi_{\cP}(S) =\bfi$ and $v_1\in V_i$ for some $i\in [s]$.
We construct absorbing $t k^2$-sets for $S$ as follows. 
We fix an edge with vertex set $W=\{v_1, x_2, \dots, x_k\}$ for $x_2, \dots, x_k \in \bigcup_{s+1\le j\le r}V_j \setminus \{y_2,\dots, y_k\}$ such that $\bfi_{\cP}(W)=\bfi_{\cP}(S)=\bfi$ and $W\cap S=\{v_1\}$. 
Note that by Lemma~\ref{lem:partition} (4), we have at least $\mu n^{k-1} - (k-1) n_1^{k-2} >\frac{\mu}2 n^{k-1}$ choices for $W$ (and $x_2, \dots, x_k$ are in $\bigcup_{s+1\le i\le r}V_i$, by the definition of $I_{\cP,2}^{s,\mu}(H)$). 
Without loss of generality, we may assume that for all $i\in \{2, \dots, k\}$, $x_i, y_i$ are in the same part $V_j$ of $\cP$, $j>s$. 
Since $x_i$ is $(\beta, t)$-reachable to $y_i$, there are at least $\beta n_1^{t k-1}$ $(t k-1)$-sets $T_i$ in $V(H)\setminus V_0$ such that both $H[T_i\cup \{x_i\}]$ and $H[T_i\cup \{y_i\}]$ have perfect matchings. 
We pick disjoint reachable $(t k-1)$-sets in $V(H)\setminus V_0$ for each $x_i, y_i$, $i\in \{2, \dots, k\}$ greedily, while avoiding the previously chosen vertices. 
Since the number of previously chosen vertices is at most $t k(k-1)+(k-1)$, we have at least $\frac{\beta}2 n_1^{t k-1}$ choices for such $(t k-1)$-sets in each step.
At last, let us pick a matching $M$ of size $t$ in $H$ that is vertex disjoint from the existing vertices (the purpose is to let the absorbing set contain exactly $tk^2$ vertices).
For the number of choices for $V(M)$, we can sequentially choose disjoint edges satisfying any $\mu$-robust edge vector $\bfi\in I_{\cP,1}^\mu(H)$ and infer that there are at least $\frac12\mu^t n^{tk}$ choices.

Note that each choice of $(W\setminus \{v_1\}) \cup T_2\cup \cdots \cup T_{k} \cup V(M)$ is an absorbing set for $S$. 
First, it contains a perfect matching because each $T_i\cup \{x_i\}$ for $i\in \{2, \dots, k\}$ spans $t$ vertex-disjoint edges and $M$ is a matching.
Second, $H[W\cup T_1\cup \cdots \cup T_{k}\cup S]$ also contains a perfect matching as each $T_i\cup \{y_i\}$ for $i\in \{2, \dots, k\}$ spans $t$ vertex-disjoint edges, $W$ is an edge and $M$ is a matching. 
There are at least $\frac{\mu}2 n_1^{k-1}$ choices for $W$ and at least $\frac{\beta}2 n_1^{t k-1}$ choices for each $T_i$ and $\frac12\mu^t n^{tk}$ choices for $V(M)$.
Thus we find at least 
\[
\frac{\mu}2 n^{k-1} \times \left(\frac{\beta}{2} n_1^{t k-1}\right)^{k-1} \times \frac12\mu^t n^{tk} \times \frac{1}{(tk^2)!} \geq \mu^{t+1} \beta^{k} n_1^{t k^2}
\] absorbing $t k^2$-sets for $S$, with vertices from $\bigcup_{s+1\le i\le r}V_i$ only.
%
%
\end{proof}

Continuing the proof of Lemma~\ref{lem:abs}, we apply Lemma~\ref{lem:derandomization} to the graph $G$ with parts $U = \{S \subseteq V(H):\bfi_{\cP}(S)\in I_{\cP}^{s, \mu}(H)\}$ and $W = \{T\subseteq \bigcup_{s+1\le i\le r} V_i: |T| = tk^2\}$, where $T_1, T_2$ in $W$ are adjacent if and only if $T_1\cap T_2 \neq \emptyset$, and $S\in U$ and $T\in W$ are adjacent if and only if $H[T]$ and $H[T\cup S]$ contain perfect matchings. In the notation of Lemma~\ref{lem:derandomization} we have $N = \binom {n_1} {tk^2}$, $M\le \binom {n_1} k + (n-n_1)\binom {n_1} {k-1}\le 2\binom {n_1} k$ and \[m = |E(G[W])| \le \binom {n_1} {tk^2} \times tk^2 \times \binom {n_1}{tk^2-1} = \frac{t^2k^4N^2}{n_1-tk^2+1}.\]
We let $\beta' = \mu^{t+1}\beta^{k+1}$, $\tau = \beta'/3$ and $r^* = (\beta')^2 n_1$.
Then by Claim~\ref{clm:abs}, for every $u\in U$, we get $\deg_G(u)\ge \beta' N $ and
\[\exp \left(\frac{\tau^2 r^*}{3 \beta'}\right) = \exp \left(\frac{(\beta')^3 n_1}{27}\right) \ge 16\binom {n_1} k \ge 8M,
\]
as $n_1$ is large enough. 
Thus by Lemma~\ref{lem:derandomization}, in time $O(N^4+MN^3) = O(n_1^{4tk^2}+n_1^k n_1^{3tk^2})= O(n^{4tk^2})$, we can find a set $R\subseteq W$ which is independent in $G$ with $(1-\tau)r^*\le |R|\le r^*$ and $|N_G(u)\cap R|\ge (\beta'-\tau-\nu)r^*$ for all $u\in U$, where
\[
\nu = \frac{2mr^*}{N^2}\le \frac{2t^2k^4 N^2 (\beta')^2 n_1}{(n_1-tk^2+1)N^2}
\le \frac{2t^2k^4 (\beta')^2 n_1}{n_1-tk^2+1} < \frac{\beta'}{3}.
\]
Note that the vertices of $R$ correspond to disjoint $tk^2$-sets of $H[\bigcup_{s+1\le i\le r} V_i]$ by the definition of $G[W]$. We now remove $tk^2$-sets in $R$ that do not have a perfect matching, and denote the resulting family of $tk^2$-sets by $\E_{abs}$. Thus, $|\E_{abs}|\le \beta n_1$, and each member of $\E_{abs}$ has a perfect matching. Moreover, every $k$-vertex set $S$ with $\bfi_{\cP}(S)\in I_{\cP}^{s, \mu}(H)$ has at least 
\[\frac{\beta'r^*}{3}\ge \frac{(\beta')^3 n_1}{3} \ge \alpha n_1
\]
absorbing $tk^2$-sets in $\E_{abs}$, as $\alpha \ll \beta, \mu$.
%
%
%
\end{proof}


\subsection{Absorption lemmas for Theorem~\ref{thm:alm_mat}}

We next give two absorption lemmas, Lemmas~\ref{lem:abs1} and~\ref{lem:S-aborb}, which will be used in the algorithmic proof of Theorem~\ref{thm:alm_mat}. 
Indeed, their non-algorithmic versions were proven in~\cite{CGHW21} in support of the non-algorithmic version of Theorem~\ref{thm:alm_mat} (\cite[Theorem 1.6]{CGHW21}).
The two absorption lemmas work for different ranges of $\ell$ and complement each other.

\begin{lemma}\label{lem:abs1}
	For $1\le \ell< \lceil 2k/3\rceil$ and $\gamma>0$, suppose $1/n\ll\alpha\ll \gamma,1/k$.
Let $H$ be an $n$-vertex $k$-graph with $\delta_\ell (H)\ge (\c+\gamma)\binom{n-\ell}{k-\ell}$.
Then in $O(n^{4k+1})$ time we can find a matching $M$ in $H$ of size $|M|\le \gamma n/k$ such that for any subset $R\subseteq V(H)\setminus V(M)$ with $|R|\le \alpha^{2} n$, $H[R\cup V(M)]$ contains a matching covering all but at most $k+1$ vertices which can be constructed in time $O(n)$.
\end{lemma}

The proof of Lemma~\ref{lem:abs1} resembles other ones by the lattice-based absorption method, so we postpone it to the appendix.

Next we give the second absorbing lemma, originally proven in~\cite{CGHW21} and we shall provide an algorithmic proof. 
Given a set $S$ of $2k-\ell$ vertices, an edge $e\in E(H)$ that is disjoint from $S$ is called \emph{$S$-absorbing} if there are two disjoint edges $e_1$ and $e_2$ in $E(H)$ such that $|e_1\cap S| = k-\lfloor \ell/2\rfloor$, $|e_1\cap e| = \lfloor \ell/2\rfloor$, $|e_2\cap S| = k-\lceil \ell/2\rceil$, and $|e_2\cap e| = \lceil \ell/2\rceil$. Consider a matching $M$ and a $(2k-\ell)$-set $S$, $V(M)\cap S = \emptyset$. If $M$ contains an $S$-absorbing edge $e$, then one can ``absorb" $S$ into $M$ by swapping $e$ for $e_1$ and $e_2$ ($k-\ell$ vertices of $e$ become uncovered). 


\begin{lemma}\label{lem:S-aborb}
	For $\lceil 2k/3 \rceil \le \ell \le k-1$ and $\gamma'>0$, suppose $1/n \ll \beta \ll \gamma', 1/k$. Let $H$ be an $n$-vertex $k$-graph with $\delta_\ell (H)\ge \gamma' n^{k-\ell}$. Then in time $O(n^{5k-\ell})$ we can find a matching $M'$ in $H$ of size $|M'|\le \beta n/k$ and such that for any subset $R\subseteq V(H)\setminus V(M')$ with $|R| \le \beta^2 n$, $H[V(M')\cup R]$ contains a matching covering all but at most $2k-\ell-1$ vertices.
\end{lemma}

\begin{proof}
		We apply Lemma~\ref{lem:derandomization} to the graph $G$ with parts $U = \binom {V(H)}{2k-\ell}$ and $W = E(H)$, where $T_1, T_2$ in $W$ are 
adjacent if and only if $T_1\cap T_2 \neq \emptyset$, and $S\in U$ and $T\in W$ are adjacent if and only if $T$ is $S$-absorbing. In the notation of Lemma~\ref{lem:derandomization} we have $N = |W| \ge \gamma' \binom {n} {k}$, $M= \binom {n} {2k-\ell}$ and 
\[m = |E(G[W])| \le N \times k \times \binom {n} {k-1}   \le \frac{k^2N^2}{\gamma'(n-k+1)}.\]
We let $\beta' = {(\gamma')^3}/{(2k!)}$, $\tau = \beta'/3$ and $r = \beta n /k$.
It is proven by simple counting in~\cite[Claim 4.1]{CGHW21} that for every $S\in \binom {V(H)}{2k-\ell}$, there are at least $ (\gamma')^3 n^k/(2k!)$ $S$-absorbing edges.
Thus in our proof, for every $u\in U$, we get $\deg_G(u)\ge \beta' n^k $ and
\[\exp \left(\frac{\tau^2 r}{3 \beta'}\right) = \exp \left(\frac{\beta' \beta n}{27 k}\right) = \exp \left(\frac{ (\gamma')^3 \beta n}{54 k! k}\right)\ge 8\binom {n} {2k-\ell} = 8M,
\]
as $n$ is large enough.
By Lemma~\ref{lem:derandomization}, in time $O(N^4+MN^3) = O(n^{4k}+n^{2k-\ell} n^{3k})= O(n^{5k-\ell})$, we can find a set $M'\subseteq W$ which is independent in $G$ with $(1-\nu)r\le |M'|\le r$ and $|N_G(u)\cap M'|\ge (\beta'-\tau-\nu)r$ for all $u\in U$, where
\[\nu = \frac{2mr}{N^2}\le \frac{2k^2 N^2 \beta n}{\gamma'(n-k+1)N^2 k }
\le \frac{2k^2 \beta n}{\gamma'(n-k+1)k } < \frac{\beta'}{3}
\]
because $\beta \ll \gamma'$.
Note that $M'$ is a matching in $H$ by the definition of $G[W]$ and  $|M'|\le \beta n/k$. 
Moreover, every $(2k-\ell)$-set $S$ has at least 
\[\frac{\beta'r}{3}\ge \frac{(\gamma')^3 \beta n}{6k! k} \ge \beta^2 n
\]
$S$-absorbing edges in $M'$, as $\beta \ll \gamma'$.

Now let $R\subseteq V(H)\setminus V(M')$ with $|R|\le \beta^2 n$ be given.
As long as there are at least $2k-\ell$ vertices uncovered, we take a set $S$ of $2k-\ell$ vertices and find an $S$-absorbing edge in $M'$.
Then we proceed the swap as described before the statement of the lemma and then reduce the number of uncovered vertices by $k$.
Applying this greedily we obtain a matching in $H[V(M')\cup R]$ covering all but at most $2k-\ell-1$ vertices.
\end{proof}


\section{An algorithmic version of Theorem~\ref{thm:alm_mat}}

In this section we give an algorithmic version of Theorem~\ref{thm:alm_mat}, which we restate as Theorem~\ref{thm:coverbut}.
In the proof, we will use the Weak Regularity Lemma, which is an extension of Szemer\'edi's regularity lemma for graphs~\cite{Sze}.

Let $H = (V, E)$ be a $k$-graph and let $A_1,\dots, A_k$ be mutually disjoint non-empty subsets of $V$. We define $e(A_1, \dots, A_k)$ to be the number of \emph{crossing edges}, namely, those with one vertex in each $A_i$, $i\in [k]$, and the density of $H$ with respect to $(A_1, \dots, A_k)$ as
\[d(A_1, \dots, A_k) = \frac{e(A_1, \dots, A_k)}{|A_1|\cdots|A_k|}.
\]
We say a $k$-tuple $(V_1, \dots, V_k)$ of mutually disjoint subsets $V_1, \dots, V_k \subseteq V$ is \emph{$(\varepsilon, d)$-regular}, for $\varepsilon > 0$ and $d\ge 0$, if 
\[ |d(A_1, \dots, A_k)-d|\le \varepsilon
\]
for all $k$-tuples of subsets $A_i \subseteq V_i$, $i\in [k]$, satisfying $|A_i|\ge \varepsilon|V_i|$. We say $(V_1, \dots, V_k)$ is \emph{$\varepsilon$-regular} if it is $(\varepsilon, d)$-regular for some $d\ge 0$. It is immediate from the definition that in an $(\varepsilon, d)$-regular $k$-tuple $(V_1, \dots, V_k)$, if $V'_i\subseteq V_i$ has size $|V'_i| \ge  c|V_i|$ for some $c\ge \varepsilon$, then $(V'_1, \dots, V'_k)$ is $(\varepsilon/c, d)$-regular.
A partition $V_0\cup V_1\cup\dots\cup V_t$ of $V$ is called \emph{$\varepsilon$-regular} if it satisfies the following properties:
\begin{enumerate}
\item $|V_0|\le \varepsilon |V|$,
\item $|V_1| = |V_2| =\cdots = |V_t|$,
\item for all but at most $\varepsilon \binom{t}{k}$ $k$-subsets $\{i_1, \dots, i_k\}\subseteq [t] $, the $k$-tuple $(V_{i_1}, \dots, V_{i_k})$ is $\varepsilon$-regular.
\end{enumerate}

We will use the following theorem, which is an algorithmic regularity lemma for weighted hypergraphs proved by Czygrinow and R\"odl~\cite{CR00}. 
For convenience, we state only the non-weighted version. 

\begin{theorem}\cite{CR00}\label{thm:regularity}
For every $k, m$, and $\varepsilon$ there exist $M, L$, and an algorithm which for any $k$-graph $H=(V, E)$ with $|V|=n\ge L$ finds in $O(n^{2k-1}\log^2 n)$ time an $\varepsilon$-regular partition $V_0\cup V_1\cup\dots\cup V_t$ of $H$ with $m\le t \le M$.
\end{theorem}

Given an $\varepsilon$-regular partition of $H$ and $d\ge 0$, we refer to $V_i, i\in [t]$ as \emph{clusters} and define the \emph{cluster hypergraph} $\mathcal K =\mathcal K(\varepsilon, d)$ with vertex set $[t]$ such that $\{i_1, \dots, i_k\}\subseteq [t]$ is an edge if and only if $(V_{i_1}, \dots, V_{i_k})$ is $\varepsilon$-regular and $d(V_{i_1}, \dots, V_{i_k})\ge d$. We will also use the following proposition.

\medskip
\begin{prop}~\cite[Proposition 16]{HS}\label{prop:cluster_hyper}
    Given an $n$-vertex $k$-graph ${H} = (V, E)$ with minimum $(k-1)$-degree
\[
\delta_{k-1}({H}) \geq \left( \frac{1}{2(k-\ell)} + \gamma \right)n
\]
and an $\varepsilon$-regular partition $V = V_0 {\cup} V_1 {\cup} \cdots {\cup} V_t$ with $0 < \varepsilon < \gamma^2/16$ and $t_0 \geq 8k/\varepsilon \geq 3k/\gamma$. Further, let $\mathcal{K} = \mathcal{K}(\varepsilon, \gamma/6)$ be the cluster hypergraph of ${H}$. Then the number of $(k-1)$-sets $S = \{i_1, \dots, i_k\} \in \binom{[t]}{k-1}$ violating
\[
\deg_{\mathcal{K}}(S) \geq \left( \frac{1}{2(k-\ell)} + \frac{\gamma}{4} \right)t
\]
is at most $\sqrt[3]{\varepsilon} \binom{t}{k}$.
\end{prop}

Combining Theorem~\ref{thm:regularity} and Proposition~\ref{prop:cluster_hyper}, we have the following corollary. It shows that the cluster hypergraph almost inherits the minimum degree of the original hypergraph. Its proof is standard and similar to the one of Proposition~\ref{prop:cluster_hyper} so we omit it.

\begin{corollary}\label{cor:cluster}
Given $c, \varepsilon, d, \ell>0$, integers $k\ge 3$ and $t_0$, there exist $T_0$ and $n_0$ such that the following holds. 
There exists an algorithm which for any $k$-graph $H=(V, E)$ on $n>n_0$ vertices with $\delta_\ell(H)\ge c\binom {n-\ell}{k-\ell}$ finds in $O(n^{2k-1}\log^2 n)$ time an $\varepsilon$-regular partition $V_0\cup V_1\cup\dots\cup V_t$ of $H$ with $t_0\le t\le T_0$. 
Moreover, in the cluster hypergraph $\mathcal K =\mathcal K(\varepsilon, d)$, all but at most $\sqrt[3] \varepsilon \binom{t}{\ell}$ $\ell$-subsets $S$ of $[t]$ satisfy $\deg_{\mathcal K} (S) \ge (c-d-\sqrt\varepsilon)\binom{t - \ell}{k -\ell}$.
\end{corollary}

We will use the following lemma to show that random subgraphs of hypergraphs typically inherit minimum-degree conditions.

\begin{lemma}\cite{FK22}\label{lem:prob}
	There is $c'=c'(k)>0$ such that the following holds. Consider an $n$-vertex $k$-graph $G$ where all but $\delta \binom n \ell$ of the $\ell$-sets have degree at least $(\mu+\eta)\binom{n-\ell}{k-\ell}$. Let $S$ be a uniformly random subset of $Q\ge 2\ell$ vertices of $G$. Then with probability at least $1-\binom Q \ell (\delta +e^{-c'\eta^2 Q})$, the random induced subgraph $G[S]$ has minimum $\ell$-degree at least $(\mu + \eta/2)\binom{Q-\ell}{k-\ell}$.
\end{lemma}

The following lemma will be used in our reproof of Theorem~\ref{thm:alm_mat}. Its non-algorithmic version appeared in~\cite{CGHW21}.
\begin{lemma}\label{lem:coverbut}
	Given integers $k, \ell$ such that $1\le \ell \le k-1$ and $\gamma, \sigma>0$, the following holds for sufficiently large $n$. Suppose $H$ is an $n$-vertex $k$-graph with $\delta_\ell(H)\ge (c_{k,\ell}^*+\gamma)\binom{n-\ell}{k-\ell}$, then in $O(n^{2k-1}\log^2 n)$ time we can find a matching $M$ in $H$ such that $M$ covers all but at most $\sigma n $ vertices.
\end{lemma}

\begin{proof}
Choose $1/n \ll 1/T_0 \ll 1/t_0\ll \varepsilon \ll 1/Q \ll \gamma, \sigma, 1/k $.
Let 
$d = \gamma/6$ and $c= c^*_{k,\ell}+\gamma$.
	By Corollary~\ref{cor:cluster}, in $O(n^{2k-1}\log^2 n)$ time, we find an $\varepsilon$-regular partition $V_0 \cup V_1\cup \dots \cup V_t$ of $H$ with $t_0\le t\le T_0$.
	 Moreover, in the cluster hypergraph $\mathcal K =\mathcal K(\varepsilon, d)$, all but at most $\sqrt[3] \varepsilon \binom{t}{\ell}$ $\ell$-subsets $S$ of $[t]$ satisfy 
	\[\deg_{\mathcal K} (S) \ge (c^*_{k,\ell} +\gamma-d-\sqrt\varepsilon)\binom{t - \ell}{k -\ell}\ge(c^*_{k,\ell} +\gamma/2)\binom{t - \ell}{k -\ell}.\]
	We randomly partition $V(\mathcal K)$ into $\lfloor t/Q \rfloor$ disjoint subsets of size $Q$ (and possibly at most $Q-1$ vertices).
	By Lemma~\ref{lem:prob}, for each subset $S$ of size $Q$, we have
	\[\Pr\left[\delta_\ell (\mathcal K[S])\le (c^*_{k,\ell} +\gamma/4)\binom{Q - \ell}{k -\ell} \right] \le \binom Q \ell (\sqrt[3] \varepsilon + e^{-c' \gamma^2Q/4}) < e^{-\sqrt Q},
	\]
	as $\varepsilon \ll 1/Q \ll \gamma, 1/k$.
	Then the expected number of subsets $S$ of size $Q$ for which $\mathcal K[S]$ has minimum $\ell$-degree at most $(c^*_{k,\ell} +\gamma/4)\binom{Q - \ell}{k -\ell}$ is no more than
	$\frac t Q e^{-\sqrt Q}$. 
	 Thus there exists a choice of $S_1, S_2, \dots, S_{t'}$ with $t'\ge (1-e^{-\sqrt Q})\frac t Q$ such that they are disjoint subsets of $V(\mathcal K)$ with $\delta_\ell (\mathcal K [S_i]) \ge (c_{k,\ell}^*+\gamma/4)\binom{Q-\ell}{k -\ell}$.
	 Moreover, in constant time we obtain $S_1, S_2, \dots, S_{t'}$ because the number of vertices of the reduced graph $\mathcal K$ is at most $T_0$ which is a constant.
	Then we can apply (the non-algorithmic version of) Theorem~\ref{thm:alm_mat} on each $\mathcal K[S_i]$, and conclude that each $\mathcal K[S_i]$ contains a matching $M_i$ in $\mathcal K$ covering $(1-\varepsilon)Q$ vertices with at most $\varepsilon Q$ vertices uncovered.
	Suppose $M_i = \{e_{i1}, e_{i2}, \dots, e_{it''}\}$ for $e_{ij} \in E(\mathcal  K)$ with $i\in [t']$ and $j\in [t'']$. Then each $e_{ij}$ corresponds to a $k$-tuple $(V_{i_1}, V_{i_2}, \dots, V_{i_k})$ in $H$ which is $\varepsilon$-regular and $d(V_{i_1}, V_{i_2}, \dots, V_{i_k})\ge d$. Then in $O(n^{k})$ time we can greedily find an almost perfect matching in the $k$-tuple $(V_{i_1}, V_{i_2}, \dots, V_{i_k})$ covering $(1-\varepsilon)k |V_{i_1}|$ vertices.
	Combining these for all edges $e_{ij}$ with $i\in [t']$ and $j\in [t'']$, we obtain a matching covering at least
	\[ (1-\varepsilon) t' t'' k |V_{i_1}| >  k (1-\varepsilon)  \times  (1-e^{-\sqrt Q}) \frac t Q  \times (1-\varepsilon) \frac Q k \times  \frac {n-|V_0|} t >  (1-\varepsilon)^3(1-e^{-\sqrt Q}) n \ge (1-\sigma) n,
	\]
	as $\varepsilon, 1/Q \ll \sigma$. In particular, the overall running time is $O(n^{2k-1}\log^2 n)$.
\end{proof}

\begin{theorem}\label{thm:coverbut}
Let $k, d$ be integers such that $1\le \ell \le k-1$ and $\gamma >0$; then there exists an $n_0\in \mathbb N$ such that the following holds for $n\ge n_0$. Suppose $H$ is an $n$-vertex $k$-graph with $\delta_\ell(H) \ge (c^*_{k, \ell} + \gamma){n-\ell \choose k-\ell}$; then in $O(n^{5k-\ell})$ time we can find a matching $M$ in $H$ that covers all but at most $2k-\ell-1$ vertices. In particular, when $n\in k\mathbb N$, $M$ is a perfect matching or covers all but exactly $k$ vertices.
\end{theorem}

For the full proof of Theorem~\ref{thm:coverbut} we refer to the proof of~\cite[Theorem 1.6]{CGHW21}.
The difference is that we use the lemmas developed in this paper (Lemmas~\ref{lem:abs1},~\ref{lem:S-aborb} and~\ref{lem:coverbut}) to replace the non-algorithmic versions in~\cite{CGHW21}. Indeed, by Lemmas~\ref{lem:abs1} and~\ref{lem:S-aborb}, in $O(n^{4k+1} +n^{5k-\ell})= O(n^{5k-\ell})$ time we can find a matching $M$ in $H$ of size $|M| \le \eta n/k$ such that, for every subset $R \subseteq  V \setminus  V (M)$ with $| R|  \leq  \sigma n$, $H[R \cup  V (M)]$ contains a matching covering all but at most $\max\{ k + 1, 2k  -  \ell  -  1\}  = 2k - \ell - 1$ vertices. Then applying Lemma~\ref{lem:coverbut} on $H_1 := H[V\setminus V(M)]$, in $O(n^{2k-1}\log^2 n)$ time we can find a matching $M_1$ in $H_1$ covering all but at most $\sigma n$ vertices of $V (H_1)$. 
Finally we can use the absorption property of $M$ to absorb the leftover vertices greedily, in time $O(n)$.


\section{Proof of Theorem~\ref{thm_structural}}

%

Now we are ready to prove Theorem~\ref{thm_structural}.
Let $H$ be an $n$-vertex $k$-graph, and let $\cP$ be a partition given by Lemma~\ref{lem:partition} satisfying (1)--(5).
We first prove the forward implication.

\subsection{Proof of the forward implication of Theorem~\ref{thm_structural}}

If $H$ contains a perfect matching $M$, then $\bfi_\cP (V(H)\setminus V(M))={\bf 0} \in L^{s, \mu}_{\cP}(H)$.
Let $M'$ be the smallest submatching of $M$ that covers $V_0$, so $|M'|\le |V_0|$.
We shall show that there exists a matching $M'' \subseteq M \setminus M'$ such that $|M''|\le q$ and $\bfi_\cP (V(H)\setminus V(M'\cup M''))\in L^{s, \mu}_{\cP}(H)$, implying that $(\cP,L^{s, \mu}_{\cP}(H))$ is $(V_0, q)$-soluble.

Indeed, suppose that $M''\subseteq M\setminus M'$ is a smallest matching such that $\bfi_{\cP}(V(H)\setminus V(M'\cup M''))\in L_{\cP}^{s, \mu}(H)$ and $|M''|=m\ge q$.
Let $M''=\{e_1, \dots, e_{m}\}$ and consider the $m+1$ partial sums
\[
\sum_{i=1}^{j}\bfi_{\cP}(e_i)+L_{\cP}^{s, \mu}(H),
\]
for $j=0,1,\dots, m$.
Since $|Q(\cP, L_{\cP}^{s, \mu}(H))|\le q\le m$, two of the sums must be in the same coset.
That is, there exist $0\le j_1< j_2\le m$ such that
\[
\sum_{i=j_1+1}^{j_2}\bfi_{\cP}(e_i) \in L_{\cP}^{s, \mu}(H).
\]
So the matching $M^*:=M''\setminus \{e_{j_1+1},\dots, e_{j_2}\}$ satisfies that $\bfi_{\cP}(V(H)\setminus V(M^*\cup M'))\in L_{\cP}^{s, \mu}(H)$ and $|M^*|< |M''|$, a contradiction.
This completes the proof of the forward implication.

\medskip
Before giving the proof of the backward implication of Theorem~\ref{thm_structural}, we first introduce the following useful constant.
Given a set $I$ of $k$-vectors in $\mathbb Z^r$, and $m^*\in \mathbb N$, consider the set $J$ of all $m'$-vectors that are in the lattice in $\mathbb Z^r$ generated by $I$ with $0\le m'\le m^*$.
That is, for any $\bfv\in J$, there exist $a_{\bfi}\in \mathbb Z$, $\bfi\in I$ such that
\[
\bfv=\sum_{\bfi\in I}a_{\bfi} \bfi.
\]
Then let $C^*:=C^*(r, k, I, m^*)$ be the maximum of $|a_{\bfi}|$ over all such $\bfv$. 
Furthermore, let 
\begin{enumerate}[label=$(\dagger)$]
\item $C_{\max}:=C_{\max}(k, m^*)$ be the maximum of $C^*=C^*(r, k, I, m^*)$ over all $r\le r_0(k):= 2^{\binom{2k-1}{k-1}}+k$ and all families of $k$-vectors $I\subseteq \mathbb Z^r$. \label{item:1}
\end{enumerate}

\subsection{Proof of the backward implication of Theorem~\ref{thm_structural}}
Recall that $c_{k,\ell}^*\ge c_{k,k-1}^* = 1/k$.
Then $\lfloor 1/c_{k,\ell}^*\rfloor \le k$.
Define constants
\[
t:=2^{k} \quad \text{ and } \quad C:= C_{\max}(k, kq+k).
\]
Define additional constants $\a, c'>0$ so that 
\[
1/n_0 \ll \alpha, 1/c' \ll \beta, \mu \ll \delta' \ll \gamma, 1/k, 1/q, 1/C, 1/t.
\]
Let $n\ge n_0$ be a multiple of $k$.
Let $H$ be as in the statement of the theorem and $\cP$ be a partition of $V(H)$ satisfying Lemma~\ref{lem:partition} (1)--(5), where the $C$ therein is as defined above.
In particular, Property~\ref{item:225} and the choice of $t$ imply that for $s+1\le i\le r$,  $|V_i|\ge \delta'n/2$ and $V_i$ is $(\beta, t)$-closed in $H[\bigcup_{s+1\le i\le r}V_i]$.
Furthermore, assume that $(\cP, L^{s, \mu}_{\cP}(H))$ is $(V_0, q)$-soluble, that is, there is a matching $M_1$ of size at most $|V_0|+q$
such that $M_1$ covers $V_0$ and it is a $(|V_0|+q)$-solution, that is,
\[
\bfi_{\cP}(V(H)\setminus V(M_1))\in L_{\cP}^{s, \mu}(H).
\]
Let $n_1:=|\bigcup_{s+1\le i\le r} V_i|$.
We first apply Lemma~\ref{lem:abs} to $H$ and get a family $\E_{abs}$ in $O(n^{4tk^2})$ time consisting of at most $\beta n_1$ disjoint $t k^2$-sets such that {$V(\E_{abs})\subseteq \bigcup_{s+1\le i\le r}V_i$ and} every $k$-set $S$ of vertices with $\bfi_{\cP}(S)\in I_{\cP}^{s, \mu}(H)$ has at least $\alpha n_1$ absorbing $t k^2$-sets in $\E_{abs}$. 

Note that $V(M_1)$ may intersect $V(\E_{abs})$ in at most $(|V_0|+q)k$ absorbing sets of $\E_{abs}$. 
Let $\E_{0}$ be the subfamily of $\E_{abs}$ obtained by removing the $t k^2$-sets that intersect $V(M_1)$. 
Let $M_0$ be the perfect matching on $V(\E_0)$ that is the union of the perfect matchings on each member of $\E_0$. 
Note that every $k$-set $S$ with $\bfi_{\cP}(S)\in I_{\cP}^{s, \mu}(H)$ has at least $\alpha n_1 - (|V_0|+q)k$ absorbing sets in $\E_{0}$.

Next we want to ``store'' some disjoint edges for each $k$-vector in $I_{\cP}^{s, \mu}(H)$ for later steps, and at the same time we also cover the remaining vertices of $\bigcup_{1\le i\le s} V_i$ (recall that $V_0$ is covered by $M_1$).
More precisely, set $C':=C^*(r, k, I_{\cP}^{s, \mu}(H), kq+k)\le C$.
Note that Lemma~\ref{lem:partition} (3) guarantees that for each $i\in [s]$, $V_i$ has at least $\binom{2k-2}{k-1} C$ uncovered vertices.
We construct a matching $M_2$ in $H - V(M_{0}\cup M_1)$ which consists of $C'$ disjoint edges $e$ with $\bfi_{\cP}(e)=\bfi$ for every $\bfi\in I_{\cP}^{s, \mu}(H)$.
So 
\[
|M_2|\le |I_{\cP}^{s, \mu}(H)| C' \le \binom{k+r-1}{k}C'.
\]
Note that $H$ contains at least $\mu n^k$ edges for each $\bfi\in I_{\cP,1}^{s, \mu}(H)$ and that
every vertex in $\bigcup_{1\le i\le s} V_i$ is in at least $\mu n^{k-1}$ edges with index vector $\bfi$ for some $\bfi\in I_{\cP,2}^{s, \mu}(H)$. 
Because 
\begin{equation}\label{eqn:bound}
	|V(M_0\cup M_{1}\cup M_2)|\le  t k^2 \alpha n_1 + (|V_0|+q )k + \binom{k+r-1}{k}C' k < \mu n_1 < \mu n,
\end{equation}
we can choose the desired edges in a greedy manner to build $M_2$.
Indeed, for every $i\in [s]$, the number of $\mu$-robust index vectors $\bfi$ such that $\bfi|_i=1$ is at most 
$\binom{(k-1)+(r-s)-1}{k-1}\le \binom{2k-2}{k-1}$ (by $r-s\le c =\lfloor 1/\delta \rfloor \le k$), and thus the process above needs at most $\binom{2k-2}{k-1}C$ uncovered vertices from $V_i$, which is okay by our construction\footnote{Remark. This is where we need Lemma~\ref{lem:partition}~\ref{item:223}, the lower bound of $|V_i|$, $i\in [s]$.}.
Finally, note that $M_2$ is constructed with $\mu$-robust edges in time $O(n^k)$.

Let $H':=H - V( M_0\cup M_1 \cup M_2)$ and $n':=|H'|$. 
So $n'\ge n - \mu n$ and 
\[
\delta_{\ell}(H') \ge \delta_{\ell}(H) - \mu n^{k-\ell} \ge  (c_{k,\ell}^*+\gamma /2) \binom{n'-\ell}{k-\ell}.
\]
%
%
By~Theorem~\ref{thm:coverbut}, in time $O(n^{5k-\ell})$ we construct a matching $M_3$ in $H'$ covering all but at most $k$ vertices.
Let $U$ be the set of vertices in $H'$ uncovered by $M_3$.
We are done if $U=\emptyset$.
Otherwise because $k$ divides $n$ we have $|U|=k$.

We write $Q:=Q(\cP, L_{\cP}^{s, \mu}(H))$ for brevity.
Recall that $\bfi_{\cP}(V(H)\setminus V(M_1))\in L_{\cP}^{s, \mu}(H)$. 
Note that by definition, the index vectors of all edges in $M_2$ are in $I_{\cP}^{s, \mu}(H)$. 
So we have $\bfi_{\cP}(V(H)\setminus V( M_1\cup M_2))\in L_\cP^{s, \mu}(H)$, namely, $R_Q(V(H)\setminus V( M_1\cup M_2))=\mathbf{0}+ L_{\cP}^{s, \mu}(H)$. 
Thus,
\[
 \sum_{e \in M_{0}\cup M_3} R_Q(e)  +R_Q(U) = \mathbf{0}+ L_{\cP}^{s, \mu}(H).
\]
Suppose $R_Q(U)=\bfv_{0}+L_{\cP}^{s, \mu}(H)$ for some $\bfv_0\in L_{\max}^{r}$; so
\[
\sum_{e \in M_{0}\cup M_3} R_Q(e)  = -\bfv_0+ L_{\cP}^{s, \mu}(H).
\]

We use the following claim proved in~\cite{HT2020} (earlier versions appeared in~\cite{Han14_poly, KKM2015}).
Its proof is via the coset arguments and is very similar to the one used in the proof of the forward implication.

\begin{claim}\label{clm:35}\cite[Claim 5.1]{HT2020}
In $O(n)$ time we can find  $e_1, \dots, e_{p} \in M_{0}\cup M_3$ for some $p\le q-1$ such that 
\begin{equation}\label{eq:RQ}
\sum_{i\in [p]} R_Q(e_i)  = -\bfv_0+ L_{\cP}^{s, \mu}(H).
\end{equation}
\end{claim}

That is, we have $\sum_{i\in [p]} \bfi_{\cP}(e_i)+\bfi_{\cP}(U)\in L_{\cP}^{s, \mu}(H)$. Let $Y:=\bigcup_{i\in [p]} e_i\cup U$ and thus $|Y|=  p k +k \leq q k +k$.
We now complete the perfect matching by absorption. Since $\bfi_{\cP}(Y)\in L_{\cP}^{s, \mu}(H)$, we have the following equation
\[
\bfi_{\cP}(Y) = \sum_{\bfv\in I_{\cP}^{s, \mu}(H)} a_{\bfv} \bfv,
\]
where $a_{\bfv}\in \mathbb{Z}$ for all $\bfv\in I_{\cP}^{s, \mu}(H)$. 
Since $|Y|\le qk+k$, by the definition of $C'$, we have $|a_{\bfv}|\le C'$ for all $\bfv\in I_{\cP}^{s, \mu}(H)$.
Noticing that $a_{\bfv}$ may be negative, we can assume $a_{\bfv}=b_{\bfv} - c_{\bfv}$ such that one of $b_{\bfv}, c_{\bfv}$ is $|a_{\bfv}|$ and the other is zero for all $\bfv\in I_{\cP}^{s, \mu}(H)$. So we have
\[
\sum_{\bfv\in I_{\cP}^{s, \mu}(H)} c_{\bfv} \bfv + \bfi_{\cP}(Y) = \sum_{\bfv\in I_{\cP}^{s, \mu}(H)} b_{\bfv} \bfv.
\]
This equation means that given a family $\E = \{W_1^{\bfv},\dots, W_{c_{\bfv}}^{\bfv}: \bfv\in I_{\cP}^{s, \mu}(H)\}$ of disjoint $k$-subsets of $V(H)\setminus Y$ such that $\bfi_{\cP}(W_{i}^{\bfv})=\bfv$ for all $i\in [c_{\bfv}]$, we can regard $V(\E)\cup Y$ as the union of disjoint $k$-sets $\{S_1^{\bfv},\dots, S_{b_{\bfv}}^{\bfv}: \bfv\in I_{\cP}^{s, \mu}(H)\}$ such that  $\bfi_{\cP}(S_j^{\bfv})=\bfv$, $j\in [b_{\bfv}]$ for all $\bfv\in I_{\cP}^{s, \mu}(H)$.
Since $c_{\bfv}\le C'$ for all $\bfv$ and $V(M_2)\cap Y=\emptyset$, we can choose the family $\E$ as a subset of $M_2$.  
In summary, starting with the matching $M_0\cup M_1\cup M_2\cup M_3$ leaving $U$ uncovered, we delete the edges $e_1, \dots, e_{p}$ from $M_{0}\cup M_3$ given by Claim~\ref{clm:35} and then leave $Y=\bigcup_{i\in [p]}V(e_i)\cup U$ uncovered.
Next we delete the family $\E$ of edges from $M_2$ and leave $V(\E)\cup Y$ uncovered.
Note that the family $\E$ can be found in constant time.
Finally, we regard $V(\E)\cup Y$ as the union of at most 
\[
|M_2|+qk+k\le \alpha n_1/2
\] 
$k$-sets $S$ each with $\bfi_{\cP}(S)\in I_{\cP}^{s, \mu}(H)$. 

Note that by definition, $Y$ may intersect at most $qk+k$ absorbing sets in $\E_{0}$, which cannot be used to absorb those sets we obtained above.
Since each $k$-set $S$ has at least $\alpha n_1- (|V_0|+q)k>\alpha n_1/2 + qk+k$ absorbing $tk^2$-sets in $\E_{0}$, we can greedily match each $S$ with a distinct absorbing $t k^2$-set $E_S\in \E_{0}$ for $S$ in time $O(n)$.
Replacing the matching on $V(E_S)$ in $M_0$ by the perfect matching on $H[E_S\cup S]$ for each $S$ gives a perfect matching in $H$.

In conclusion, our algorithm outputs a perfect matching of $H$ in time $O(n^{4tk^2})=O(n^{2^{k+2}k^2})$.

\section*{Acknowledgements}
The second author would like to thank Andrew Treglown and Hi{\d{\^{e}}}p H\`an for helpful discussions at the early stage of this project.
The authors would like to thank two anonymous referees for their helpful comments that improved the presentation of this paper.

\bibliographystyle{abbrv}
\bibliography{Bibref, refs}

\appendix
\section{Proof of Lemma~\ref{lem:abs1}}

We first present two propositions useful in the proof of Lemma~\ref{lem:abs1}. For $i\in[r]$, let $\bold{u}_i \in \mathbb{Z}^r$ be the $i$-th \emph{unit vector}, namely, $\bold{u}_i$ has $1$ on the $i$-th coordinate and $0$ on other coordinates. A \emph{transferral} is the vector $\bold{u}_i-\bold{u}_j$ for some $i\neq j$.
To merge some pair of parts, we need to use~\cite[Lemma 3.4]{Han15_mat} to keep the closedness of $V_i$ and $V_j$ if $\bold{u}_i-\bold{u}_j\in L_{\cP}^{0,\mu}(H)$.
In particular, we can test whether the edge-lattice contains any transferral in constant time.

\begin{proposition}\label{prop:absorb}\cite{CGHW21}
	Given $\min\{3, k/2\} \leq \ell< \lceil 2k/3 \rceil$, suppose $1/n\ll\mu \ll \varepsilon \ll \gamma$.
Let $H$ be an $n$-vertex $k$-graph with $\delta_{\ell}(H) \ge(1/4+\gamma)\binom{n-\ell}{k-\ell}$, and let $\cP = \{V_0,V_1, \dots, V_r \}$ be a partition of $V(H)$ with $r \leq 3$ such that $|V_{0}|\leq \sqrt{\varepsilon}n$ and for each $i \in [r], |V_{i}|\geq \varepsilon^{2}n$, and $L_{\cP}^{0, \mu}(H)$ contains no transferral.
Then for every $U\subseteq V(H)\setminus V_{0}$ with $|U|=k+2$, there exist $i,j \in [r]$ such that $\bfi_{\cP}(U)-\bfu_{i}-\bfu_{j}\in L_{\cP}^{0, \mu}(H)$.
\end{proposition} 

\begin{proposition}\label{prop:52}\cite[Proposition 3.6]{Han15_mat}
 Given $\min\{3,k/2\}\le \ell< k-2$ or $(k, \ell) = (5, 2)$, suppose that $1/n\ll \mu \ll \varepsilon \ll \gamma$.
 Let $H$ be an $n$-vertex $k$-graph with $\delta_\ell (H) \ge (1/3+\gamma)\binom {n-\ell}{k-\ell}$ and let $\cP = \{V_0, \dots, V_r\}$ be a partition of $V(H)$ with $r\le 3$ such that $|V_0|\le \sqrt{\varepsilon} n$ and for any $i\in [r]$, $|V_i|\ge \varepsilon^2 n$, and $L_{\cP}^{0, \mu} (H)$ contains no transferral. Then for any $U\subseteq V(H)\setminus V_0$ with $|U| = k+1$, there exists $i\in[r]$ such that $\bfi_{\cP} (U) -\bfu_i \in L_{\cP}^{0, \mu} (H)$.
\end{proposition}

Note that both propositions above guarantee the existence of certain vector in the robust lattice, which can be found in constant time. We also need the following result.


\begin{lemma}\cite[Lemma 3.4]{Han15_mat}\label{lem:transferral}
	Given $0<\mu, \beta \ll \varepsilon \ll 1/i_0$, there exist $0<\beta'\ll \mu, \beta$ and an integer $t\ge i_0$ such that the following holds for sufficiently large $n$. Suppose $H$ is an $n$-vertex $k$-graph, and $\cP = \{V_0, V_1, \dots, V_r\}$ is a partition with $r\le i_0$ such that $|V_0|\le \sqrt{\varepsilon}n$ and for any $i\in[r]$, $|V_i|\ge \varepsilon^2 n$ and $V_i$ is $(\beta, i_0)$-closed in $H$. If {$\bold{u}_i-\bold{u}_j\in L_{\cP}^{0,\mu}(H)$} for $i, j \in [r]$ and $i\neq j$, then $V_i\cup V_j$ is $(\beta', t)$-closed in $H$.
\end{lemma}

Now we present the proof of Lemma~\ref{lem:abs1}.

\begin{proof}[Proof of Lemma~\ref{lem:abs1}]
We split into three cases according to the values of $k$ and $\ell$.

\noindent{\bf{Case I:}} $\ell=1$ and $k\ge 3$, or $\ell=2$ and $k\ge 6$. Let $\alpha:=\gamma/(k-1)!$. Notice that we have
$c_{k,1}^*\ge 1-(1-1/k)^{k-1} \ge 1-(2/3)^2 > 1/2$ and $c_{k,2}^*\ge 1-(1-1/k)^{k-2} \ge 1-(5/6)^4 > 1/2$.
That is, in this case we always have $\delta_\ell (H)\ge (\frac 1 2 + \gamma)\binom {n-\ell}{k-\ell}$. 
Note that every two vertices are $(\alpha, 1)$-reachable in $H$ because for $u, v\in H$, $|N(u)|\ge (1/2+\gamma){n-1 \choose k-1}$ and $|N(v)|\ge (1/2+\gamma){n-1 \choose k-1}$ such that $|N(u)\cap N(v)|\ge \gamma {n-1 \choose k-1} =\frac {\gamma}{(k-1)!} n^{k-1}$. Thus $H$ is $(\alpha, 1)$-closed.
We apply Lemma~\ref{lem:abs} on $H$ with the trivial partition $\cP=\{V_1\}$ to obtain a family $\mathcal{F}_{abs}$ of disjoint $k^2$-sets. 
In particular, we have $|V(\mathcal{F}_{abs})|\le k^2\beta n$, $H[V(\mathcal{F}_{abs})]$ contains a perfect matching, and every $k$-vertex set $S $ has at least $\alpha n$ absorbing $k^2$-set in $\mathcal{F}_{abs}$.
For any subset $R\subseteq V(H)\setminus V(\mathcal{F}_{abs})$, $H[R\cup \mathcal{F}_{abs}]$ has a perfect matching and the matching can be constructed in $O(n)$ time by absorbing the vertices in $R$.

\medskip
\noindent{\bf{Case II:}} $\min\{3, k/2\} \leq \ell< \lceil 2k/3 \rceil$.
Note that we have $\c\ge 1-(1-1/k)^{k/3}\ge 1-e^{-1/3}>1/4$.
Given $\gamma>0$. Let $C:=C_{\max}(k,2k)$, where $C_{\max}(k,m)$ is defined in~\ref{item:1}.
We first define additional constants as
\[
1/n\ll \alpha\ll \beta'\ll\beta,\mu\ll \varepsilon \ll \gamma,1/k,1/t,1/C.
\]
Let $H$ be a $k$-graph with $\delta_\ell(H)\ge (\c+\gamma){n-\ell \choose k-\ell}$. Since $\ell>1$, by Proposition~\ref{prop_deg}, we have $\delta_1(H)\ge (1/4+\gamma)\binom{n-1}{k-1}$ from which we conclude that among every four vertices, there exist two vertices that are $(\alpha,1)$-reachable. 
Applying Lemma~\ref{S-closed} with $c=3$, in time $O(n^{k+1})$ we can find a set of vertices $S\subseteq V(H)$ with $|S|\geq (1-3\delta' )n$ such that $|\tilde{N}_{\a, 1}(v, H[S])|\geq \delta' n$ for any $v\in S$. Let $V_0 : = V(H)\setminus S$.
Then applying Lemma~\ref{lem:P} on $S$ with $\delta'=\c+\gamma/2>1/4+\gamma/2$, in $O(n^{4 k+1})$ time we can find a partition $\mathcal{P}'=\{V'_1,\dots,V'_{r'}\}$ of $S$ with $r' \le \min\{3, 1/\delta'\} = 3$. In particular, for any $i\in[r']$, $|V'_i|\ge \varepsilon^2 n$ and $V'_i$ is $(\beta,8)$-closed in $H$. 
Let $\cP = \{V_0, V'_1,\dots,V'_{r'}\}$.
Next we consider transferrals in the $\mu$-robust edge-lattice.
We greedily check whether $\mathbf{u}_i-\mathbf{u}_j\in L_{\mathcal{P}}^{0, \mu}(H)$ for distinct $i, j\in[r']$. 
If so, then we merge $V'_i$ and $V'_j$, $i, j>0$ to one part.
By Lemma~\ref{lem:transferral}, $V'_i\cup V'_j$ is $(\beta'', t')$-closed for some $\beta''> 0$ and $t' \ge 4$.
Note that in $O(n^k)$ time, we greedily check the existence of transferrals and merge the parts until there is no transferral in the $\mu$-robust edge-lattice.
Let $\mathcal{P}_0 = \{V_0,V_{1},\dots, V_{r}\}$ be the resulting partition for some $1\le r\le 3$.
Note that we may have applied Lemma~\ref{lem:transferral} at most twice, and we see that $V_i$ is $(\beta', t)$-closed for each $i\in [r]$ by the choice of $\beta'$.
Then we apply Lemma~\ref{lem:abs} on $H$ with partition $\cP_0$ with $s=0$ and obtain a family $\mathcal{F}_{abs}$ of disjoint $tk^2$-sets satisfying that $|V(\mathcal{F}_{abs})|\le tk^2\beta n$, $H[V(\mathcal{F}_{abs})]$ contains a perfect matching $M_1$, and every $k$-vertex set $S$ with $\mathbf{i}_{\mathcal{P}_0}(S)\in I^{0, \mu}_{\mathcal{P}_0}(H)$ has at least $\alpha(n-|V_0|)$ absorbing $tk^2$-set in $\mathcal{F}_{abs}$.

Let $V'=V\setminus V(\mathcal{F}_{abs})$ and $H'=H[V']$. Now we find a matching $M_2$ in $H'$ as follows.
For each $\mathbf{v} \in I^{0, \mu}_{\mathcal{P}_0}(H)$, we greedily pick a matching $M_{\bfv}$ of size $C\alpha^2n$ such that $\mathbf{i}_{\mathcal{P}_0}(e) = \mathbf{v}$ for every $e\in M_{\bfv}$.
Then let $M_2$ be the union of $M_{\bfv}$ for all $\mathbf{v}\in I^{0, \mu}_{\mathcal{P}_0}(H)$, and we have $V_0\cap V(M_2)=\emptyset$.
It is possible to pick $M_2$ because there are at least $\mu n^k$ edges $e$ with $\mathbf{i}_{\mathcal{P}_0}(e) = \mathbf{v} \in I^{0, \mu}_{\mathcal{P}_0}(H)$.
To be more precise, since $|I^{0, \mu}_{\mathcal{P}_0}(H)|\le \binom{k+r-s-1}{r-s-1}\le \binom{k+2}{2}$ and $\alpha\ll \beta' \ll \mu \ll 1/t,1/C$, we obtain
\[
|V(M_2)\cup V(\mathcal{F}_{abs})| \le k|I^{0, \mu}_{\mathcal{P}_0}(H)|C\alpha^2 n+tk^2 \beta n<\mu n,
\]
which yields that the number of edges intersecting these vertices is less than $\mu n^k$, as required.

Next we build a matching $M_3$ to cover all vertices in $V_0\setminus V(\mathcal{F}_{abs})$.
Note that $|M_3|\le |V_0|\le \sqrt{\varepsilon}n$.
Specifically, when we greedily match a vertex $v\in V_0\setminus V(\mathcal{F}_{abs})$, we need to avoid at most $k|M_3|+|V(M_2)\cup V(\mathcal{F}_{abs})|\le k\sqrt{\varepsilon}n+\mu n\le 2k\sqrt{\varepsilon}n$ vertices, and thus at most $2k\sqrt{\varepsilon}n^{k-1}$ $(k-1)$-sets.
Since $\delta_1(H)> \gamma\binom{n-1}{k-1}>2k\sqrt{\varepsilon}n^{k-1}$, we can always find a desired edge containing $v$ and put it to $M_3$ as needed.
Since both $M_2$ and $M_3$ are constructed greedily, the running time is $O(n^k)$.

Let $M=M_1\cup M_2\cup M_3$.
It is easy to see that $M$ is a matching because $M_1, M_2$ and $M_3$ are pairwise vertex disjoint.
Now we prove that $M$ is a desired matching satisfying the conclusion of Lemma~\ref{lem:abs1}.
Note that $|M|\leq|M_3|+|V(M_2)\cup V(\mathcal{F}_{abs})|/k\le 2\sqrt{\varepsilon}n\leq \gamma n/k$ since $\varepsilon\ll \gamma,1/k$.
Consider any subset $R\subseteq V\setminus V(M)$ with $|R|\le \alpha^2 n$.
We are already done if $|R|\le k+1$.
Otherwise fix any set $U\subseteq R$ of $k+2$ vertices, there exist $i,j\in [r]$ such that $\mathbf{i}_{\mathcal{P}_0}(U)- \mathbf{u}_i-\mathbf{u}_j\in L^{0, \mu}_{\mathcal{P}_0}(H)$ by Proposition~\ref{prop:absorb}.
Note that this does not guarantee that we can delete one vertex $u$ from $U\cap V_i$ and delete another vertex $v$ from $U\cap V_j$ such that $\mathbf{i}_{\mathcal{P}_0}(U\setminus \{u,v\})\in L^{0, \mu}_{\mathcal{P}_0}(H)$, because it is possible that $U \cap V_i=\emptyset$ or $U \cap V_j=\emptyset$ for the $i,j$ returned by the proposition.
As $\ell\ge 2$, there is a vector $\mathbf{v}\in I^{0, \mu}_{\mathcal{P}_0}(H)$ such that $\mathbf{v}_{V_i} \ge 1$ and $\mathbf{v}_{V_j} \ge 1$.
Notice that $M_2$ contains $C\alpha^2 n$ edges with index vector $\mathbf{v}$.
Fix one such edge $e\in E_{\mathbf{v}}$ and two vertices $v_1\in e\cap V_i$, $v_2\in e\cap V_j$.
We delete $e$ from $M_2$ and let $U' = U\cup (e\setminus \{v_1,v_2\})$.
Clearly, $\mathbf{i}_{\mathcal{P}_0}(U')\in L^{0, \mu}_{\mathcal{P}_0}(H)$ and $|U'| = 2k$.
Hence, by the definition of $L^{0, \mu}_{\mathcal{P}_0}(H)$, there exist nonnegative integers $b_{\mathbf{v}}, c_{\mathbf{v}}$ for all $\mathbf{v}\in I^{0, \mu}_{\mathcal{P}_0}(H)$ such that
\[
\mathbf{i}_{\mathcal{P}_0}(U')=\sum\limits_{\mathbf{v}\in I^{0, \mu}_{\mathcal{P}_0}(H)}b_{\mathbf{v}} \mathbf{v}-\sum\limits_{\mathbf{v}\in I^{0, \mu}_{\mathcal{P}_0}(H)}c_{\mathbf{v}} \mathbf{v},\]
which implies that
\[
\mathbf{i}_{\mathcal{P}_0}(U')+\sum\limits_{\mathbf{v}\in I^{0, \mu}_{\mathcal{P}_0}(H)}c_{\mathbf{v}} \mathbf{v}=\sum\limits_{\mathbf{v}\in I^{0, \mu}_{\mathcal{P}_0}(H)}b_{\mathbf{v}} \mathbf{v}.
\]
We have that $b_{\mathbf{v}}, c_{\mathbf{v}}\le C$ from the definition of $C$.
For each $\mathbf{v}\in I^{0, \mu}_{\mathcal{P}_0}(H)$, we pick $c_{\mathbf{v}}$ edges in $M_2$ with index vector $\mathbf{v}$.
By the equation above, the union of these edges and $U'$ can be partitioned as a collection of $k$-sets, which contains exactly $b_{\mathbf{v}}$ $k$-sets $F$ with $\mathbf{i}_{\mathcal{P}_0}(F) = \mathbf{v}$ for each $\mathbf{v} \in I^{0, \mu}_{\mathcal{P}_0}(H)$.
We repeat the process at most $\alpha^2 n/k$ times until there are at most $k+1$ vertices left.
Note that for each $\mathbf{v}\in I^{0,\mu}_{\mathcal{P}_0}(H)$, our algorithm consumes at most $(1 + C)\alpha^2n/k < C\alpha^2n$ edges from $M_2$ with index vector $\mathbf{v}$, which is possible by the definition of $M_2$.
Furthermore, after the process, we obtain at most $\left(2+C|I^{0, \mu}_{\mathcal{P}_0}(H)|\right)\alpha^2n/k\le \left(2+C\binom{k+2}{2}\right)\alpha^2 n/k<\alpha n$ $k$-sets $S$ with $\mathbf{i}_{\mathcal{P}_0}(S)\in I^{0, \mu}_{\mathcal{P}_0}(H)$ since $\alpha\ll 1/k, 1/C$.
By the absorbing property of $\mathcal{F}_{abs}$, we can greedily absorb them by $\mathcal{F}_{abs}$ and get a matching $M_4$.
Thus, $H[R\cup V(M)]$ contains a matching covering all but at most $k + 1$ vertices.

\medskip
\noindent{\bf{Case III:}} $\ell =2$ and $k=5$. Notice that $c_{5,2}^* = 1-(1-1/5)^3 >1/3$. The proof is identical to Case II except that we replace Proposition~\ref{prop:absorb} with Proposition~\ref{prop:52}.
\end{proof}

\end{document}